\documentclass[EJP,preprint]{ejpecp} 

\DeclareMathAlphabet{\mathbbold}{U}{bbold}{m}{n}
\usepackage{mathrsfs}
\usepackage{xfrac}
\allowdisplaybreaks[1]

\SHORTTITLE{Talagrand's Bound for Product Semigroups}

\TITLE{Talagrand’s Convolution Bound for Tensor Powers of Finite-State Markov Semigroups}

\AUTHORS{%
  Alexander~Marsh\footnote{Colorado State University, United States of America.
    \EMAIL{marshajk@colostate.edu}}\orcid{0009-0002-9235-6161}}


\AMSSUBJ{60E15 ; 60J27 ; 60J35}

\ABSTRACT{Talagrand's convolution conjecture, recently proven
by J. Lu, S. Guo, and E. Fang, concerns the smoothing properties of the heat semigroup on $\{-1,1\}^N$. We follow their perturbed coupled reverse-heat approach to generalize the theorem to tensor powers of irreducible continuous-time Markov semigroups on finite state spaces. Specifically, let $(S_t)_{t\geq0}$ be the semigroup of an irreducible continuous-time Markov chain on a finite state space $G_0$ with stationary distribution $\pi$. Then, for every $\tau>0$, there exists a finite constant $C_{\tau,S,\pi}<\infty$ such that, for every $N\geq1$, every nonnegative function $f:G_0^N\to[0,\infty)$, and every $\eta>1$,
\begin{align*}
    \pi^{\otimes N}\left(S_\tau^{\otimes N}f>\eta\int_{G_0^N}f\,d\pi^{\otimes N}\right)\leq\frac{C_{\tau,S,\pi}}{\eta\sqrt{\log\eta}},
\end{align*}
where the constant is independent of $N$, $f$, and $\eta$. Most of the machinery developed by Lu, Guo, and Fang transfers cleanly, but since probability distributions on $|G_0|$ elements are not determined by their mean for $|G_0|>2$, we must use explicit terminal bridge laws to find a dimension-independent bound for the terminal-discrepancy estimate. We have formalized our result in Lean~4. The formalization is publicly available.\footnote{\url{https://github.com/AMarsh00/TalagrandConvolutionExtension}}}

\begin{document}
\raggedbottom



\section{Introduction}

Markov's inequality gives a general tail bound of
\begin{align*}
    \pi^{\otimes N}\left(f>\eta\int_{G}f\,d\pi^{\otimes N}\right)\leq \frac{1}{\eta},\qquad \eta>1,
\end{align*}
for every nonnegative function $f$ on $G$, where $G:=G_0^N$ for a finite state space $G_0$, $N\geq 1$, and a fully-supported probability measure $\pi$ on $G_0$. A natural question, posed originally by Talagrand \cite{Talagrand1989Conjecture} for the $\sfrac{\mathbb{Z}}{2\mathbb{Z}}$ case, is whether applying a tensor power of a Markov semigroup with stationary distribution $\pi$ to $f$ improves this bound. Following Chen's remark~\cite[page 2]{chen2025talagrand}, it is trivial to do this if we allow the constant to depend on $N$. Our goal is to find a constant that does not depend on $N$. This would show that tensor powers of semigroups have very nice smoothing properties regardless of the power $N$. The $\sfrac{\mathbb{Z}}{2\mathbb{Z}}$ case was known as Talagrand's convolution conjecture~\cite{Talagrand1989Conjecture}, positing that at each fixed time, the $1/\eta$ Markov bound can be improved to $C/(\eta\sqrt{\log\eta})$, for $C$ independent of $N$. Lu, Guo, and Fang~\cite{lu2026weaktype} recently established this bound, building on earlier work on reverse heat processes and perturbed couplings by Chen~\cite{chen2025talagrand}, and Xiang and Zhang \cite{xiang2026layerwise}. A subsequent martingale proof of the Boolean bound was recently discovered by Shaposhnikov~\cite{shaposhnikov2026remarks}.

In this paper, we extend Lu, Guo, and Fang's result to tensor powers of any irreducible continuous-time Markov semigroup with any stationary distribution with full support on a finite state space of any size. For each fixed time, we obtain the same $1/(\eta\sqrt{\log\eta})$ tail decay, with a constant depending on the semigroup and stationary distribution but independent of $N$. We first establish the bound for a specific semigroup which we call the $\pi$-refresh semigroup. We then obtain the general result by noticing that we can factor any semigroup through the $\pi$-refresh semigroup. Explicitly, our main theorem is as follows:

\begin{theorem}\label{thm:1}
    Let $(S_t)_{t\geq0}$ be the semigroup of an irreducible
    continuous-time Markov chain on $G_0$, with stationary distribution $\pi$. Then, for every $\tau>0$, there exists a constant $C_{\tau ,S,\pi}<\infty$ such that, for every $N\geq1$, every nonnegative $f:G\to[0,\infty)$, and every $\eta>1$,
    \begin{align*}
        \pi^{\otimes N}\left(S_\tau^{\otimes N}f>          \eta\int_{G}f\,d\pi^{\otimes N}\right)\leq\frac{C_{\tau ,S,\pi}}{\eta\sqrt{\log\eta}}.
    \end{align*}
    In particular, the constant is independent of $N$.
\end{theorem}

Our proof follows the approach used by Lu, Guo, and Fang~\cite{lu2026weaktype}, which itself uses the perturbed reverse coupling ideas of Chen~\cite{chen2025talagrand} and the localized terminal discrepancy method of Xiang and Zhang~\cite{xiang2026layerwise}. The main challenge in the extension is that the Boolean proof relies on the fact that a probability distribution is determined by its mean, which is clearly no longer true for state spaces of size larger than $2$. To fix this, we instead compute the full conditional terminal laws and bound their changes under a coordinate replacement using a variance inequality which removes the $N$ dependence. The other parts of our argument are very similar to those previously presented by Chen, Lu, Guo, and Fang, and are re-proven here as our extended coupling construction makes it somewhat unclear that they will transfer without viewing the details.

We have formally verified the main theorem in Lean 4, allowing some larger constants in our intermediate lemmas. We represent our coupled jump processes in the same way that the Lean formalization accompanying Lu, Guo, and Fang~\cite{lu2026weaktype} does. The source code and verification instructions are available in the accompanying repository, which is linked in a footnote on our first page.

\section{Main proof}
To prove Theorem~\ref{thm:1}, we first show the result for a specific semigroup that we call the $\pi$-refresh semigroup. This semigroup is a conceptual analogue of the Boolean heat semigroup on more general finite state spaces. Given $t\geq0$ and $x\in G$, let $X_t^x$ be a random vector on $G$ whose coordinates are independent and satisfy
\begin{equation*}
    X_t^{x,[i]}\sim e^{-t}\delta_{x_i}+(1-e^{-t})\pi.
\end{equation*}
Equivalently, independently in each coordinate, $X_t^{x,[i]}$ remains equal to $x_i$ with probability $e^{-t}$ and is replaced by an independent sample from $\pi$ with probability $1-e^{-t}$. Thus, for every $a\in G_0$,
\begin{equation*}
    \mathbb{P}\left[X_t^{x,[i]}=a\right]=e^{-t}\mathbbold{1}_{\{a=x_i\}}+(1-e^{-t})\pi(a).
\end{equation*}
The $\pi$-refresh semigroup $(P_t^\pi)_{t\geq0}$ is then defined, for every function $f:G\to[0,\infty)$ and every $x\in G$, by
\begin{equation}\label{pi_refresh}
    P_t^\pi f(x):=\mathbb{E}\left[f(X_t^x)\right]=\int_G f(y)\ d\pi_{t,x}(y),
\end{equation}
where $\pi_{t,x}:=\bigotimes\limits_{i=1}^N\left(e^{-t}\delta_{x_i}+(1-e^{-t})\pi\right)$, using the notation that $\delta_{x_i}$ is the Dirac delta measure that assigns full probability to $x_i$ and no probability to other elements. The proof of the following lemma is deferred to Section~\ref{sec3}.

\begin{lemma}[$\pi$-refresh tail bound]\label{lem:lem1}
    For all $\tau>0$, there exists a constant $C_{\tau ,\pi}<\infty$ depending only on $\tau$ and $\pi$ such that, for every $N\geq 1$, every nonnegative $f:G\to[0,\infty)$, and every $\eta>1$,
    \begin{align*}
        \pi^{\otimes N}\left(P_\tau^\pi f>\eta\int_{G}f\ d\pi^{\otimes N}\right)\leq\frac{C_{\tau,\pi}}{\eta\sqrt{\log\eta}}.
    \end{align*}
\end{lemma}

From Lemma~\ref{lem:lem1}, we can prove our main theorem.

\begin{proof}[Proof of Theorem~\ref{thm:1}]
    Let $K_\tau$ be the transition matrix of $S_\tau$ at time $\tau$, i.e. $K_\tau(x,y)=S_\tau\mathbbold{1}_y(x)$, $x,y\in G_0$. Since $\pi$ is the stationary distribution of an irreducible Markov chain, $\pi(x)>0$ for all $x\in G_0$. Let $\epsilon_\tau:=\min\limits_{x,y\in G_0}\frac{K_\tau(x,y)}{\pi(y)}$. Since the underlying continuous-time Markov chain is irreducible on a finite state space, $K_\tau(x,y)>0$ for every $\tau>0$ and all $x,y\in G_0$. Hence $\epsilon_\tau>0$. Moreover, $K_\tau(x,y)\geq\epsilon_\tau\pi(y)$, so summing over $y$ gives $1\geq\epsilon_\tau$. $\epsilon_\tau=1$ if and only if $K_\tau(x,y)=\pi(y)$ for all $x,y\in G_0$, in which case $S_\tau^{\otimes N}f(x)=\int_G f\ d\pi^{\otimes N}$ for all $x$, so for $\eta>1$ the event $S_\tau^{\otimes N}f>\eta\int_G f\ d\pi^{\otimes N}$ never occurs. Thus assume $0<\epsilon_\tau<1$ and define $R_\tau(x,y):=\frac{K_\tau(x,y)-\epsilon_\tau\pi(y)}{1-\epsilon_\tau}$, noting that $R_\tau$ is stochastic and thus a valid transition matrix. We then let $T_\tau$ be the Markov operator of $R_\tau$. The stationary distribution of $T_\tau$ is $\pi$ since by stationarity of $K_\tau$,
    \begin{align*}
        \sum\limits_{x\in G_0}\pi(x) R_\tau(x,y)=\sum\limits_{x\in G_0}\frac{\pi(x) K_\tau(x,y)-\epsilon_\tau\pi(x) \pi(y)}{1-\epsilon_\tau}=\frac{\pi(y)-\epsilon_\tau\pi(y)}{1-\epsilon_\tau}=\pi(y).
    \end{align*}
    Let $\Pi f:=\int_{G_0}f\ d\pi$. Then $(1-\epsilon_\tau)T_\tau f+\epsilon_\tau\Pi f=S_\tau f$ for all $f:G_0\to[0,\infty)$. By stationarity, $\Pi=\Pi T_\tau$. Using this and tensoring, we have that for all $f:G\to[0,\infty)$,
    \begin{align*}
        ((1-\epsilon_\tau)I+\epsilon_\tau\Pi)^{\otimes N}T_\tau^{\otimes N} f=S_\tau^{\otimes N}f.
    \end{align*}
    Now using the definition of $P_t$~\eqref{pi_refresh} and independence of the coordinates,
    \begin{align*}
        P_t^\pi
        =
        \bigl(e^{-t}I+(1-e^{-t})\Pi\bigr)^{\otimes N}.
    \end{align*}
    If we set $t=-\log(1-\epsilon_\tau)$, this matches $(1-\epsilon_\tau)I+\epsilon_\tau\Pi$. Therefore,
    \begin{align*}
        S_\tau^{\otimes N}f=P_t^\pi T_\tau^{\otimes N}f.
    \end{align*}
    We then apply Lemma \ref{lem:lem1} to $T_\tau^{\otimes N}f$ (noting that $T_\tau^{\otimes N}f\geq0$ because $T_\tau$ is a Markov operator) to find
    \begin{align*}
        \pi^{\otimes N}\left(P_t^\pi T_\tau^{\otimes N}f>\eta\int_GT_\tau^{\otimes N}f\ d\pi^{\otimes N}\right)\leq\frac{C_{t ,\pi}}{\eta\sqrt{\log\eta}}.
    \end{align*}
    Finally, as $\pi$ is stationary for $T_\tau$, the $T_\tau$ in the integral can be removed, leaving only
    \begin{align*}
        \pi^{\otimes N}\left(S_\tau^{\otimes N}f>\eta\int_Gf\ d\pi^{\otimes N}\right)\leq\frac{C_{\tau,S ,\pi}}{\eta\sqrt{\log\eta}},
    \end{align*}
    where $C_{\tau,S,\pi}:=C_{-\log(1-\epsilon_\tau),\pi}$.
\end{proof}

\begin{remark}
    We note that the proof above does not require $(S_t)_{t\geq0}$ to be a semigroup, so this result holds for any strictly positive finite stochastic kernel with stationary distribution $\pi$.
\end{remark}

\section{$\pi$-refresh semigroup tail bound}\label{sec3}
To prove Lemma~\ref{lem:lem1}, we follow the same strategy as Lu, Guo, and Fang~\cite{lu2026weaktype}. For any $f:G\to(0,\infty)$ with $\int_G f\ d\pi^{\otimes N}=1$, define $\nu_f:=f\cdot\pi^{\otimes N}$ (i.e. $f$ is the Radon-Nikodym derivative $f=\frac{d\nu_f}{d\pi^{\otimes N}}$). We use $\lesssim_\pi$ for an inequality up to a constant depending only on $\pi$. We use the same definition of the anti-concentration profile as Chen~\cite[Equation 11]{chen2025talagrand} for a function $f:G\to(0,\infty)$ such that $\int_G f\ d\pi^{\otimes N}=1$. This is (for intervals $I\subseteq\mathbb{R}$)
\begin{equation}
    \mathfrak{A}_t^f(I):=\mathbb{P}_{Y\sim \nu_{P_t^\pi f}}[\log P_t^\pi f(Y)\in I].
\end{equation}
The analogue of ~\cite[Proposition 3.1]{lu2026weaktype} is
\begin{lemma}\label{lem:lem2}
    Let $\tau>0$, $f:G\to(0,\infty)$, $\int_Gf\ d\pi^{\otimes N}=1$, and $\ell>2$. Then
    \begin{align*}
        \mathfrak{A}_\tau^f((\ell,\ell+1])\lesssim_\pi\kappa_{\tau,\pi}^2\left(\frac{\log\kappa_{\tau,\pi}}{\kappa_{\tau,\pi}-1}\right)^{\frac{1}{2}}\frac{1}{\sqrt{\ell}},
    \end{align*}
    where $\kappa_{\tau,\pi}:=\frac{e^{-\tau}+(1-e^{-\tau})\max\limits_{y\in G_0}\pi(y)}{(1-e^{-\tau})\min\limits_{y\in G_0}\pi(y)}$.
\end{lemma}

This lemma is proven at the end of this section. It allows us to directly prove Lemma~\ref{lem:lem1}. Our argument follows the strategy of the proof of Theorem 1.1 in ~\cite{lu2026weaktype}.

\begin{proof}[Proof of Lemma~\ref{lem:lem1}]
    First, suppose that $f>0$ and $\int_G f\ d\pi^{\otimes N}=1$. Then the goal is to show that $\pi^{\otimes N}(P_\tau^\pi f>\eta)\leq \frac{C_{\tau,\pi}}{\eta\sqrt{\log\eta}}$. For $\eta>e^2$, we expand by using Lemma~\ref{lem:lem2} (letting $A_{\tau,\pi}$ be the constant in the lemma):
    \begin{align*}
        \pi^{\otimes N}(P_\tau^\pi f>\eta)&=\sum\limits_{k=0}^\infty\pi^{\otimes N}(P_\tau^\pi f\in(e^k\eta,e^{k+1}\eta])\leq\sum\limits_{k=0}^\infty\frac{1}{e^k\eta}\mathbb{E}^{\pi^{\otimes N}}\left[P_\tau^\pi f\mathbbold{1}_{P_\tau f\in(e^k\eta,e^{k+1}\eta]}\right]\\
        &=\sum\limits_{k=0}^\infty\frac{1}{e^k\eta}\mathfrak{A}_\tau^f((\log(e^k\eta),\log(e^k\eta)+1])\leq A_{\tau,\pi}\sum\limits_{k=0}^\infty\frac{1}{e^k\eta\sqrt{\log(e^k\eta)}}\\
        &\leq\frac{A_{\tau,\pi}}{\eta\sqrt{\log \eta}}\sum\limits_{k=0}^\infty\frac{1}{e^k}=\frac{\frac{A_{\tau,\pi}}{1-e^{-1}}}{\eta\sqrt{\log\eta}}.
    \end{align*}
    For $1<\eta\leq e^2$, we use Markov's inequality: $\pi^{\otimes N}(P_\tau^\pi f>\eta)\leq\frac{1}{\eta}\leq\frac{\sqrt{2}}{\eta\sqrt{\log\eta}}$. The result then follows with $C_{\tau,\pi}=\max\{\sqrt{2},\frac{A_{\tau,\pi}}{1-\frac{1}{e}}\}$. For a general nonnegative $f$, put $m=\int_Gf\,d\pi^{\otimes N}$. If $m=0$, full support forces $f$ to also be $0$, and the tail event never happens. If $m>0$, apply the normalized bound to $f_\epsilon=\frac{f+\epsilon}{m+\epsilon}$. Since $\{P_\tau^\pi f_\epsilon>\eta\}=\{P_\tau^\pi f>\eta m+(\eta-1)\epsilon\}$, these events increase to $\{P_\tau^\pi f>\eta m\}$ as $\epsilon\downarrow0$, giving the result.
\end{proof}

The proof of Lemma~\ref{lem:lem2} relies on Lemmas~\ref{lem:time-smoothed}--\ref{lem:stopped-score}, which are the general versions of Lemmas 3.2--3.5 in Lu, Guo, and Fang~\cite{lu2026weaktype}. The proofs of these lemmas are deferred to Sections \ref{sec4}, \ref{sec5}, \ref{sec6}, and \ref{sec7} respectively. These lemmas rely on a coupling construction similar to the one used by Lu, Guo, and Fang, which we build below. We use notation that for $x\in G$, $g\in G_0$, $x^{i\to g}$ is the element of $G$ matching $x$ in all coordinates except coordinate $i$, where it is instead $g$. Throughout, we can assume that $|G_0|>1$ since if $G_0$ has only one element, every Markov semigroup is the identity and we find the overall theorem trivially. For the jump SDEs below, we define subtraction on $G$ by identifying it with the Abelian group $\left(\sfrac{\mathbb{Z}}{|G_0|\mathbb{Z}}\right)^N$. Subtraction lets us clearly notate a jump $x\to x^{i\to g}$ by $x^{i\to g}-x$ in the jump equations.

\begin{lemma}\label{lem:vt}
    Let $f:G\to(0,\infty)$, $\int_Gf\ d\pi^{\otimes N}=1$. For any $T>0$, let $(\Omega,\mathcal{F},\{\mathcal{F}_t\}_{t\in[0,T]},\mathbb{P})$ be a filtered probability space containing $\{N^{i,g}\}_{i\in\{1,\cdots,N\},g\in G_0}$, which are independent Poisson random measures (PRMs) on $(0,\infty)\times\mathbb{R}$ with intensity $dt\ dz$. Let $U_0\sim\nu_f$ be independent of $\{N^{i,g}\}_{i,g}$, and take $(\mathcal{F}_t)_{t\in[0,T]}$ to be the usual augmentation of the filtration generated by $U_0$ and these PRMs. The jump SDE
    \begin{align*}
        dU_t&=\sum\limits_{i=1}^N\sum\limits_{g\in G_0}\left(U_{t^-}^{i\to g}-U_{t^-}\right)\int_\mathbb{R}\mathbbold{1}_{0< z\leq\pi(g)}\ N^{i,g}(dt,dz),
    \end{align*}
    has semigroup $(P_t^\pi)_{t\in[0,T]}$. Now let $(\Omega,\widetilde{\mathcal{F}},\{\widetilde{\mathcal{F}_t}\}_{t\in[0,T]},\mathbb{P})$ be a filtered probability space containing $\{\widetilde{N}^{i,g}\}_{i\in\{1,\cdots,N\},g\in G_0}$, which are independent PRMs on $(0,\infty)\times\mathbb{R}$ with intensity $dt\ dz$. Let $V_0\sim\nu_{P_T^\pi f}$ be independent of $\{\widetilde{N}^{i,g}\}_{i,g}$, and take $(\widetilde{\mathcal{F}}_t)_{t\in[0,T]}$ to be the usual augmentation of the filtration generated by $V_0$ and these PRMs. The jump SDE
    \begin{align*}
        dV_t&=\sum\limits_{i=1}^N\sum\limits_{g\in G_0}\left(V_{t^-}^{i\to g}-V_{t^-}\right)\int_\mathbb{R}\mathbbold{1}_{0< z\leq\pi(g)\rho_{t,i,g}(V_{t^-})}\ \widetilde{N}^{i,g}(dt,dz),
    \end{align*}
    is the time-reversed version of $U_t$, where $p_t:=P_t^\pi f$ is the density of the law of $U_t$ with respect to $\pi^{\otimes N}$ and $\rho_{t,i,g}(x)=\frac{p_{T-t}(x^{i\to g})}{p_{T-t}(x)}$.
\end{lemma}
\begin{proof}
    For any $T>0$, existence and uniqueness of $V_t$ and $U_t$ follow from the bounded total jump rate for pure-jump processes~\cite[Chapter 4.7]{ethier2009markov}. Since $U_t$ is a jump process, its generator takes the form $\mathcal{L}f(x)=\sum\limits_{x'\neq x}q(x,x')(f(x')-f(x))$, where $q(x,x')$ is the rate of jumping from $x$ to $x'$. From the SDE formula, the jump rate from $x$ to $x^{i\to g}$ is (for all $i,g\neq x_i$)
    \begin{align*}
        q(x,x^{i\to g})=\int_\mathbb{R}\mathbbold{1}_{0\leq z\leq \pi(g)}\ dz=\pi(g).
    \end{align*}
    All other $N^{i,g}$ jump rates are zero. Thus the generator of $U_t$ is
    \begin{align*}
        \sum\limits_{i=1}^N\sum\limits_{g\in G_0}\pi(g)(f(x^{i\to g})-f(x)),
    \end{align*}
    so $U_t$ has semigroup $P_t^\pi$ (by Lemma~\ref{generator}). Similarly, the generator of $V_t$ is $\mathcal{L}_tf(x)=\sum\limits_{x'\neq x}q_t(x,x')(f(x')-f(x))$ where $q_t(x,x')$ is the rate of jumping from $x$ to $x'$ at time $t$ (this is dependent on $t$ since $V_t$ has jump rates depending on time). From the SDE for $V_t$, the rate of jumping from $x$ to $x^{i\to g}$ is (for all $i,g\neq x_i$)
    \begin{align*}
        q_t(x,x^{i\to g})=\int_{\mathbb{R}}\mathbbold{1}_{0\leq z\leq \pi(g)\rho_{t,i,g}(x)}\ dz=\pi(g)\rho_{t,i,g}(x)=\pi(g)\frac{p_{T-t}(x^{i\to g})}{p_{T-t}(x)}.
    \end{align*}
    All other $\widetilde{N}^{i,g}$ jump rates are again zero. Thus the generator for $V_t$ is
    \begin{align*}
        \sum\limits_{i=1}^N\sum\limits_{g\in G_0}\pi(g)\left(\frac{p_{T-t}(x^{i\to g})}{p_{T-t}(x)}\right)(f(x^{i\to g})-f(x)),
    \end{align*}
    so the generator of $V_t$ matches the generator of the time-reversed version of $P_t^\pi$ (by Lemma~\ref{reverse_generator}), meaning that $V_t$ is the time-reversed version of $U_t$.
\end{proof}

For the coupling construction below, fix $\tau>0$, $T_0>1$, set $T=T_0+\tau$, and let $\theta\in[T_0-1,T_0)$. We use the notation that $(x)_+=\max\{x,0\}$. Fix $\ell>0$, set $\alpha=5$, and define
\begin{align*}
    R_\theta&:=\left[\ell-\log P_{T-\theta}^\pi f(V_\theta)\right]_+,\qquad\mathcal E_\theta:=\{R_\theta\geq2\alpha\},\qquad\overline\delta:=\frac{\alpha\mathbbold{1}_{\mathcal E_\theta}}{R_\theta+1}.
\end{align*}
Define the stopping time
\begin{align*}
    \sigma&:=\inf\left\{t\in[\theta,T_0]:\log P_{T-t}^\pi f(V_t)\geq\ell+1\right\}\wedge T_0.
\end{align*}
We now construct a perturbed process $(W_t)_{\theta\leq t\leq T_0}$ that is coupled to $V_t$. For $\delta\in[0,1]$ and $r>0$, define the switched-power map
\begin{align*}
    \Psi_\delta(r)&:=\begin{cases}
        r+1-r^\delta,&0<r\leq1,\\[1mm]
        r^{1-\delta},&r>1.
    \end{cases}
\end{align*}
Let $\delta_t:=\overline{\delta}\,\mathbbold{1}_{\{t\leq\sigma\}}$, where $\overline{\delta}\in[0,1]$ is $\widetilde{\mathcal F}_\theta$-measurable and $\sigma$ is the stopping time defined above. We set $W_\theta=V_\theta$ and (using the same PRMs as $V_t$)
\begin{align*}
    dW_t&=\sum_{i=1}^N\sum_{g\in G_0}\left(W_{t^-}^{i\to g}-W_{t^-}\right)\int_{\mathbb{R}}\mathbbold{1}_{\left\{0<z\leq\pi(g)\Psi_{\delta_t}\left(\rho_{t,i,g}(V_{t^-})\right)\right\}}\widetilde N^{i,g}(dt,dz),\qquad t\in(\theta,T_0].
\end{align*}
This coupling construction is the natural generalization of the construction used by Lu, Guo, and Fang~\cite{lu2026weaktype} to state spaces with arbitrary size.

\begin{lemma}
    $W_t$ is well-defined and there exists a unique strong solution to the coupled SDE $(V_t,W_t)$.
\end{lemma}
\begin{proof}
    By Lemma~\ref{lem:rho_bound}, $0<\rho_{t,i,g}(x)\le\kappa_{\tau,\pi}$ for $\theta\le t\le T_0$. Moreover,
    \begin{align*}
        0\leq\bar\delta=\frac{\alpha\mathbf{1}_{\{R_\theta\ge2\alpha\}}}{R_\theta+1}\leq\frac{\alpha}{2\alpha+1}<1.
    \end{align*}
    For $d\in[0,1]$ and $r>0$, we have $0<\Psi_d(r)\le\max\{1,r\}$. To see this, if $r\leq 1$, then $r\leq r^d\leq 1$, so $r\leq r+1-r^d\leq 1$, and if $r>1$, then $r^{1-d}\leq r$. Consequently,
    \begin{align*}
        \pi(g)\Psi_{\delta_t}\bigl(\rho_{t,i,g}(V_{t-})\bigr)\leq\pi(g)\kappa_{\tau,\pi}.
    \end{align*}
    Lemma~\ref{lem:vt} gives uniqueness of the solution $V_t$. Since $\bar\delta$ is $\widetilde{\mathcal F}_\theta$-measurable and $t\mapsto\mathbf{1}_{\{t\le\sigma\}}$ is left-continuous and adapted, $\delta_t$ is predictable on $(\theta,T_0]$, so the acceptance thresholds ($z\in(0,\pi(g)\Psi_{\delta_t}(\rho_{t,i,g}(V_{t^-}))]$) defining $W_t$ are predictable. The total jump intensity is uniformly bounded by
    \begin{align*}
        \sum_{i=1}^N\sum_{g\in G_0}\pi(g)\Psi_{\delta_t}\bigl(\rho_{t,i,g}(V_{t-})\bigr)\leq N\kappa_{\tau,\pi}.
    \end{align*}
    The compensation formula then gives an adapted c\`adl\`ag strong solution with finitely many jumps.
    
    Since the update times and coordinate replacements are determined entirely by $V_t$ and the driving PRMs, any two solutions with the same initial state and driving noise agree successively at every update. This plus uniqueness of $V_t$ gives uniqueness of $(V_t,W_t)$.
\end{proof}

We now state the lemmas we will use to prove Lemma~\ref{lem:lem2}.

\begin{lemma}[Time-smoothed anti-concentration profile estimate]\label{lem:time-smoothed}
    Let $f:G\to(0,\infty)$,\\$\int_Gf\ d\pi^{\otimes N}=1$, $\ell>2$. Then
    \begin{align*}
        \int\limits_0^\infty\mathfrak{A}_s^f((\ell,\ell+1])\ ds\lesssim_\pi\frac{1}{\ell}.
    \end{align*}
\end{lemma}

\begin{lemma}[Localized total variation distance]\label{lem:localized-tv}
    For every $A\in\widetilde{\mathcal{F}}_\theta$ such that $A\subseteq\mathcal{E}_\theta$, define
    \begin{align*}
        \mathscr S_A&:=\mathbb{E}\left[\mathbbold{1}_A\int\limits_\theta^{\sigma}\bar{\delta}^2\sum_{i=1}^N\sum_{g\in G_0}\pi(g)\left(\rho_{t,i,g}(V_{t^-})-1\right)^2\ dt\right],\\
        D_A&:=\sup_{B\subseteq G}\left|\mathbb{P}\left(A\cap\{W_{T_0}\in B\}\right)-\mathbb{P}\left(A\cap\{V_{T_0}\in B\}\right)\right|.
    \end{align*}
    Then if we let $\Lambda_{\tau,\pi}:=\frac{\kappa_{\tau,\pi}\log\kappa_{\tau,\pi}}{\kappa_{\tau,\pi}-1}$,
    \begin{align*}
        D_A\lesssim_\pi     \kappa_{\tau,\pi}\Lambda_{\tau,\pi}\sqrt{\mathscr{S}_A\mathbb{P}(A)}+\kappa_{\tau,\pi}\Lambda_{\tau,\pi}^2\mathscr{S}_A.
    \end{align*}
\end{lemma}

\begin{lemma}[Strong monotonicity]
\label{lem:strong-monotonicity}
    For $r\in\mathbb{R}$, define 
    \begin{align*}
        A_\theta(r):=\mathbb{P}\left(\mathcal{E}_\theta\cap\log P_\tau^\pi f(V_{T_0})\in(r,r+1]\right).
    \end{align*}
    Let $c_0:=e^{1-\alpha}=e^{-4}$. Then
    \begin{align*}
        (1-c_0)A_\theta(\ell)\leq c_0\sum_{j\geq1}e^{-j}A_\theta(\ell+j)+2D_{\mathcal{E}_\theta}.
    \end{align*}
\end{lemma}

\begin{lemma}[Stopped score energy]
\label{lem:stopped-score}
    We have
    \begin{align*}
        &\mathbb{E}\left[\left.\int\limits_\theta^\sigma\sum_{i=1}^N\sum_{g\in G_0}\pi(g)\left(\rho_{t,i,g}(V_{t^-})-1\right)^2\,dt\right|\widetilde{\mathcal F}_\theta\right]\leq\frac{2(\kappa_{\tau,\pi}-1)}{\log\kappa_{\tau,\pi}}\left(R_\theta+1+\log\kappa_{\tau,\pi}\right).
    \end{align*}
    Moreover, for every
    $A\in\widetilde{\mathcal F}_\theta$ satisfying
    $A\subseteq\mathcal E_\theta$,
    \begin{align*}
        \mathscr S_A&\lesssim\frac{\kappa_{\tau,\pi}}{\Lambda_{\tau,\pi}}\mathbb{E}\left[\frac{\mathbbold{1}_A}{R_\theta+1}\right]+(\kappa_{\tau,\pi}-1)\mathbb{E}\left[\frac{\mathbbold{1}_A}{(R_\theta+1)^2}\right].
    \end{align*}
\end{lemma}

We now use the proof of Lu, Guo, and Fang~\cite[Proposition 3.1]{lu2026weaktype} to prove Lemma~\ref{lem:lem2}. We write out the proof here for completeness.

\begin{proof}[Proof of Lemma~\ref{lem:lem2}]
    Let $K:=\sqrt{\kappa_{\tau,\pi}^3\Lambda_{\tau,\pi}}$. The inequalities $\frac{\kappa_{\tau,\pi}-1}{\kappa_{\tau,\pi}}\leq\log\kappa_{\tau,\pi}\leq\kappa_{\tau,\pi}-1$ give $1\leq\Lambda_{\tau,\pi}\leq\kappa_{\tau,\pi}$ and $K\geq1$. Assume first that $\ell\geq64K^2$.
    
    Set $m:=\lfloor\frac{\ell}{2}\rfloor$ and partition $\mathcal E_\theta$ into bands
    \begin{align*}
        \mathcal E_r(\theta):=\{r\leq R_\theta<r+1\},\quad 2\alpha\leq r<m,\qquad\mathcal E_\infty(\theta):=\{R_\theta\geq m\},
    \end{align*}
    $r\in\mathbb{Z}$. Lemmas~\ref{lem:localized-tv} and~\ref{lem:stopped-score} imply
    \begin{align*}
        D_{\mathcal E_r(\theta)}\lesssim_\pi\left(\frac{K}{\sqrt{r+1}}+\frac{\kappa_{\tau,\pi}^2\Lambda_{\tau,\pi}}{r+1}+\frac{\kappa_{\tau,\pi}^2\Lambda_{\tau,\pi}^2}{(r+1)^2}\right)\mathbb P[\mathcal E_r(\theta)],
    \end{align*}
    and $D_{\mathcal E_\infty(\theta)}\lesssim_\pi\frac{K}{\sqrt\ell}+\frac{\kappa_{\tau,\pi}^2\Lambda_{\tau,\pi}}{\ell}+\frac{\kappa_{\tau,\pi}^2\Lambda_{\tau,\pi}^2}{\ell^2}$. By the marginal law of $V_\theta$ and
    Lemma~\ref{lem:time-smoothed},
    \begin{align*}
        \int\limits_{T_0-1}^{T_0}\mathbb P[\mathcal E_r(\theta)]\,d\theta=\int\limits_\tau^{\tau+1}\mathfrak A_s^f((\ell-r-1,\ell-r])\,ds\lesssim_\pi\frac1\ell.
    \end{align*}
    Summing over the layers then gives (using $K^2=\kappa_{\tau,\pi}^3\Lambda_{\tau,\pi}\leq\kappa_{\tau,\pi}^4$ and monotonicity of $\frac{\log(ex)}{\sqrt{x}}$)
    \begin{equation}\label{eq:averaged-active-discrepancy}
        \int\limits_{T_0-1}^{T_0}D_{\mathcal E_\theta}\,d\theta\lesssim_\pi\frac{K}{\sqrt\ell}+\frac{\kappa_{\tau,\pi}^2\Lambda_{\tau,\pi}\log(e\ell)}{\ell}+\frac{\kappa_{\tau,\pi}^2\Lambda_{\tau,\pi}^2}{\ell}\lesssim_\pi\frac{K}{\sqrt\ell}.
    \end{equation}
    We now define
    \begin{align*}
        \overline{A_\theta}(\ell):=\mathbb P\left[\mathcal E_\theta^c,\log P_\tau^\pi f(V_{T_0})\in(\ell,\ell+1]\right].
    \end{align*}
    On $\mathcal{E}_\theta^c$, by definition we have $\log P_{T-\theta}^\pi f(V_\theta)>\ell-2\alpha$. Using Lemma~\ref{lem:time-smoothed}, if the log density starts in $(\ell-2\alpha,\ell+2]$, the $\overline{A_\theta}$ term contributes $\mathcal{O}_\pi(\frac{1}{\ell})$ after integrating from $T_0-1$ to $T_0$. $M_t:=\frac{1}{P_{T-t}^\pi f(V_t)}$ is a bounded martingale since its time derivative and reverse-generator drift cancel. Conditional Markov's inequality therefore gives
    \begin{align*}
        \mathbb P\left[\log P_\tau^\pi f(V_{T_0})\leq\ell+1\ \middle|\ \widetilde{\mathcal F}_\theta\right]
        \leq\frac{e^{\ell+1}}{P_{T-\theta}^\pi f(V_\theta)}.
    \end{align*}
    On the starting band $(\ell+2+k,\ell+3+k]$, the $\overline{A_\theta}$ term is at most $e^{-1-k}$. Applying Lemma~\ref{lem:time-smoothed} to these bands gives
    \begin{equation}\label{eq:averaged-inactive-probability}
        \int\limits_{T_0-1}^{T_0}\overline A_\theta(\ell)\,d\theta\lesssim_\pi\frac1\ell+\sum_{k\geq0}\frac{e^{-1-k}}{\ell+2+k}\lesssim_\pi\frac1\ell.
    \end{equation}
    Let $B(u):=\mathfrak A_\tau^f((u,u+1])$. Since $B(\ell)=A_\theta(\ell)+\overline A_\theta(\ell)$ and $A_\theta(\ell+j)\leq B(\ell+j)$, averaging Lemma~\ref{lem:strong-monotonicity} and using \eqref{eq:averaged-active-discrepancy} and \eqref{eq:averaged-inactive-probability} yields
    \begin{align*}
        B(\ell)\leq C_\pi\frac{K}{\sqrt\ell}+\frac{c_0}{1-c_0}\sum_{j\geq1}e^{-j}B(\ell+j),\qquad \ell\geq64K^2,
    \end{align*}
    where $c_0=e^{-4}$ and $C_\pi$ depends only on $\pi$. For fixed $\ell\geq64K^2$, the quantity $\mathcal M_\ell:=\sup\limits_{u\geq\ell}\sqrt u\,B(u)$ is finite because $G$ is finite. Applying the recurrence at every $u\geq\ell$ gives
    \begin{align*}
        \mathcal M_\ell\leq C_\pi K+\frac{e^{-4}}{(1-e^{-4})(e-1)}\mathcal M_\ell.
    \end{align*}
    The coefficient of $\mathcal M_\ell$ is less than one, so $B(\ell)\lesssim_\pi \frac{K}{\sqrt{\ell}}$. For $2<\ell<64K^2$, the same bound follows from $B(\ell)\leq1\leq\frac{8K}{\sqrt\ell}$.
\end{proof}

\section{Time-smoothed anti-concentration profile estimate}\label{sec4}
Here we prove Lemma~\ref{lem:time-smoothed}. Our proof follows Chen's proof~\cite[Lemma 4]{chen2025talagrand} using a more general log-sobolev inequality.

\begin{proof}[Proof of Lemma~\ref{lem:time-smoothed}]
    Write $\mu:=\pi^{\otimes N}$ and $p_s:=P_s^\pi f$. For $s>0$, let
    \begin{align*}
        E:=\{x\in G:\log p_s(x)\in(\ell,\ell+1]\},\qquad h_s(x):=p_s(x)\chi(\log p_s(x)-\ell)^2,
    \end{align*}
    where $\chi:\mathbb{R}\to[0,1]$ is a $C^1$ cutoff satisfying $\chi=1$ on $[0,1]$, $\operatorname{supp}(\chi)\subseteq(-1,2)$, and $\|\chi'\|_\infty\leq c'$ for a universal constant $c'$. Then $0\leq h_s\leq p_s$, $h_s=p_s$ on $E$, and $h_s=0$ outside $\{\ell-1<\log p_s<\ell+2\}$. We match Chen's entropy definition
    \begin{align*}
        \operatorname{Ent}_\mu(h)=\sup\left\{\int_G h\varphi\,d\mu\ \middle|\ \int_G e^\varphi\,d\mu\leq1\right\}.
    \end{align*}
    Taking $\varphi=\ell\mathbbold{1}_E-\log\bigl(1+(e^\ell-1)\mu(E)\bigr)$, we have $\int_G e^\varphi\,d\mu=\frac{e^\ell\mu(E)+1-\mu(E)}{1+(e^\ell-1)\mu(E)}=1$. Since $\int_G h_s\,d\mu\leq\int_G p_s\,d\mu=1$,
    \begin{align*}
        \operatorname{Ent}_\mu(h_s)\geq\ell\int_E p_s\,d\mu-\log\bigl(1+(e^\ell-1)\mu(E)\bigr).
    \end{align*}
    Moreover, $\int_E p_s\,d\mu=\mathfrak{A}_s^f((\ell,\ell+1])$ and $\mu(E)\leq e^{-\ell}\mathfrak{A}_s^f((\ell,\ell+1])$. This along with $\log(1+x)\leq x$ gives
    \begin{equation}\label{eq55}
        \operatorname{Ent}_\mu(h_s)\geq(\ell-1)\mathfrak{A}_s^f((\ell,\ell+1])\geq\frac{\ell}{2}\mathfrak{A}_s^f((\ell,\ell+1]).
    \end{equation}
    Let $m_\pi:=\min\limits_{x\in G_0}\pi(x)$. The general log-Sobolev inequality of Diaconis and Saloff-Coste~\cite[Theorem A.1]{diaconis1996logarithmic}, together with tensorization, gives
    \begin{align*}
        \operatorname{Ent}_\mu(v^2)\leq c_\pi\mathcal{E}(v,v),
        \qquad c_\pi:=\frac{\log(1-m_\pi)-\log(m_\pi)}{1-2m_\pi},
    \end{align*}
    with the convention that $c_\pi=2$ when $m_\pi=\frac{1}{2}$ and  $\mathcal{E}(v,w):=-\int_G vLw\,d\mu$. Applying this inequality to $v=\sqrt{h_s}$ and using Equation~\ref{eq55} yields
    \begin{equation}\label{eq66}
        \mathfrak{A}_s^f((\ell,\ell+1])\leq\frac{2c_\pi}{\ell}\mathcal{E}(\sqrt{h_s},\sqrt{h_s}).
    \end{equation}
    To bound the Dirichlet form, set $\psi(u):=e^{\frac{u}{2}}\chi(u-\ell)$, noting that $\sqrt{h_s}=\psi(\log p_s)$. We have
    \begin{align*}
        |\psi'(u)|\leq\left(\frac12+c'\right)e^{\frac{u}{2}}\mathbbold{1}_{(\ell-1,\ell+2)}(u).
    \end{align*}
    Consequently, by Cauchy--Schwarz, for $a<b$,
    \begin{align*}
        |\psi(b)-\psi(a)|^2&\lesssim|e^b-e^a|\int\limits_{\ell-1}^{\ell+2}\mathbbold{1}_{[a,b)}(u)\,du.
    \end{align*}
    The corresponding bound for $b<a$ follows by symmetry. By reversibility of $L$,
    \begin{align*}
        &\mathcal{E}(\mathbbold{1}_{p_s>e^u},p_s)=\frac{1}{2}\sum_{i=1}^N\sum_{g\in G_0}\pi(g)\int_G\left|\mathbbold{1}_{p_s(x^{i\to g})>e^u}-\mathbbold{1}_{p_s(x)>e^u}\right||p_s(x^{i\to g})-p_s(x)|\,d\mu(x).
    \end{align*}
    Applying the preceding bound for $\psi$ with $a=\log p_s(x)$ and $b=\log p_s(x^{i\to g})$ and summing over $i\in\{1,\cdots,N\},g\in G_0$ gives a bound for the Dirichlet form, which we can use with Equation~\ref{eq66} to get
    \begin{equation}\label{eq77}
        \mathfrak{A}_s^f((\ell,\ell+1])\lesssim_\pi\frac{1}{\ell}\int\limits_{\ell-1}^{\ell+2}\mathcal{E}(\mathbbold{1}_{p_s>e^u},p_s)\,du.
    \end{equation}
    Finally, define $F_s(u):=\int_G(p_s-e^u)_+\,d\mu$. For each fixed $u$, this function is locally absolutely continuous in $s$, and the chain rule gives, for almost every $s$,
    \begin{align*}
        \frac{d}{ds}F_s(u)=\int_G\mathbbold{1}_{p_s>e^u}Lp_s\,d\mu=-\mathcal{E}(\mathbbold{1}_{p_s>e^u},p_s).
    \end{align*}
    Since $p_s\to1$ as $s\to\infty$ and $u>\ell-1>0$, we have $F_\infty(u)=0$. Also, $F_0(u)\leq\int_G f\,d\mu=1$. Finally, integrating Equation~\ref{eq77} and using Tonelli's theorem to swap integrals, we conclude that
    \begin{align*}
        \int\limits_0^\infty\mathfrak{A}_s^f((\ell,\ell+1])\,ds&\lesssim_\pi\frac{1}{\ell}\int\limits_{\ell-1}^{\ell+2}\bigl(F_0(u)-F_\infty(u)\bigr)\,du\lesssim_\pi \frac1\ell.
    \end{align*}
\end{proof}

\section{Local total variation bound}\label{sec5}
In this section we prove Lemma~\ref{lem:localized-tv} using the terminal discrepancy method of Xiang and Zhang~\cite[Section 3]{xiang2026layerwise}, adapted to our construction using similar ideas to those used to resolve the Boolean case by Lu–Guo–Fang~\cite[Section 4.1]{lu2026weaktype}. The main additional step is the full bridge estimate in Lemma~\ref{lem:coupling-energies}, which does not follow from the Boolean proof.

\begin{lemma}\label{lem:synchronous-conditioned}
    Let $\varphi:G\to[0,1]$ and fix $\zeta\in G$. Set $\mathbb P^\zeta=\mathbb P(\,\cdot\mid V_T=\zeta)$ and $H_t^\zeta(x)=\mathbb P(V_T=\zeta\mid V_t=x)$. For $0\leq t\leq T_0$ and $x,y\in G$, let $\widehat{\mathbb P}_{t,x,y}^\zeta$ denote the law of $(\widehat V_s,\widehat W_s)_{s\in[t,T_0]}$ started from $(x,y)$ at time $t$, with generator
    \begin{align*}
        \mathcal L_s^\zeta F(v,w)&=\sum_{i=1}^N\sum_{g\in G_0}\pi(g)\lambda_{s,i,g}^\zeta(v)[F(v^{i\to g},w^{i\to g})-F(v,w)],
    \end{align*}
    where $\lambda_{s,i,g}^\zeta(v)=\rho_{s,i,g}(v)\frac{H_s^\zeta(v^{i\to g})}{H_s^\zeta(v)}$. Its first coordinate has the law of $(V_s)_{s\in[t,T_0]}$ under $\mathbb P^\zeta(\,\cdot\mid V_t=x)$. Using $\widehat{\mathbb{E}}_{t,x,y}^\zeta$ to mean expectation under the law $\widehat{\mathbb{P}}_{t,x,y}^\zeta$, define
    \begin{align*}
        q_t^\zeta(x,y):=\widehat{\mathbb E}_{t,x,y}^\zeta[\varphi(\widehat W_{T_0})].
    \end{align*}
    Then $q_{T_0}^\zeta(x,y)=\varphi(y)$, $(\partial_t+\mathcal L_t^\zeta)q_t^\zeta=0$ for $t<T_0$, and $q_t^\zeta(x,x)=\mathbb E^\zeta[\varphi(V_{T_0})\mid V_t=x]$.
\end{lemma}
\begin{proof}
    Clearly $H_s^\zeta>0$ for $s\in[0,T_0]$. Since the state space is finite and $T-T_0=\tau>0$, the rates are continuous and bounded on this interval. This means that a unique solution exists for $(\widehat{V}_t,\widehat{W}_t)$. Conditioning $V_t$ on $V_T=\zeta$ multiplies its rate from $v$ to $v^{i\to g}$ by $\frac{H_s^\zeta(v^{i\to g})}{H_s^\zeta(v)}$, which can be seen from an application of Bayes' Theorem. Hence its conditioned generator matches the generator above, and uniqueness gives that the laws of $V_t$ and $\widehat{V}_t$ match. The terminal condition follows from the definition of $q$, and the finite-state backward Kolmogorov equation gives $(\partial_t+\mathcal L_t^\zeta)q_t^\zeta=0$. Finally, since the updates to $\widehat{V}_t$ and $\widehat{W}_t$ are synchronous, if $\widehat{V}_t=\widehat{W}_t$, that equality is preserved for all $s\geq t$, which proves the last identity.
\end{proof}

\begin{lemma}\label{lem:conditioned-joint-generator}
    Fix $\zeta\in G$ such that $\mathbb{P}[V_T=\zeta]>0$. For every bounded function $\varphi:G\times G\to\mathbb R$, the process
    \begin{align*}
        M_t^{\zeta,\varphi}&:=\varphi(V_t,W_t)-\varphi(V_\theta,W_\theta)-\int_\theta^t\mathcal G_s^\zeta\varphi(V_{s^-},W_{s^-})\,ds
    \end{align*}
    is a $\mathbb P^\zeta$-martingale on $[\theta,T_0]$, where
    \begin{align*}
        \mathcal G_t^\zeta\varphi(x,y)=\sum_{i=1}^N\sum_{g\in G_0}\Big[&r_{t,i,g}^\zeta(x)\bigl(a_{t,i,g}(x)\wedge b_{t,i,g}(x)\bigr)\Delta_{i,g}^{xy}\varphi(x,y)
        \\&+r_{t,i,g}^\zeta(x)\bigl(a_{t,i,g}(x)-b_{t,i,g}(x)\bigr)_+\Delta_{i,g}^{x}\varphi(x,y)
        \\&+\bigl(b_{t,i,g}(x)-a_{t,i,g}(x)\bigr)_+\Delta_{i,g}^{y}\varphi(x,y)\Big].
    \end{align*}
    Here, $a_{t,i,g}(x):=\pi(g)\rho_{t,i,g}(x),b_{t,i,g}(x):=\pi(g)\Psi_{\delta_t}\left(\rho_{t,i,g}(x)\right)$.The difference operators are
    \begin{align*}
        \Delta_{i,g}^{xy}\varphi(x,y)&:=\varphi(x^{i\to g},y^{i\to g})-\varphi(x,y),\\
        \Delta_{i,g}^{x}\varphi(x,y)&:=\varphi(x^{i\to g},y)-\varphi(x,y),\\
        \Delta_{i,g}^{y}\varphi(x,y)&:=\varphi(x,y^{i\to g})-\varphi(x,y).
    \end{align*}
    Finally, $r_{t,i,g}^\zeta(x):=\frac{H_t^\zeta(x^{i\to g})}{H_t^\zeta(x)}$. Then $\mathcal G_t^\zeta$ is the predictable joint generator of $(V_t,W_t)$ under $\mathbb P^\zeta$.
\end{lemma}
\begin{proof}
    Our first goal is to find the Radon-Nikodym derivative $\frac{d\mathbb{P}^\zeta}{d\mathbb{P}}|_{\widetilde{\mathcal{F}}_t}$. At time $t=T$, this is a function $g:\Omega\to[0,\infty)$ such that
    \begin{align*}
        \mathbb{P}^\zeta[A]=\int_A g\ d\mathbb{P}=\mathbb{E}^\mathbb{P}[\mathbbold{1}_Ag].
    \end{align*}
    Using Bayes' Theorem, we find that
    \begin{align*}
        \mathbb{P}^\zeta[A]=\mathbb{P}[A\mid V_T=\zeta]=\frac{\mathbb{P}[A\cap V_T=\zeta]}{\mathbb{P}[V_T=\zeta]}=\frac{\mathbb{E}^\mathbb{P}[\mathbbold{1}_A\mathbbold{1}_{V_T=\zeta}]}{\mathbb{P}[V_T=\zeta]}=\mathbb{E}^\mathbb{P}\left[\mathbbold{1}_A\frac{\mathbbold{1}_{V_T=\zeta}}{\mathbb{P}[V_T=\zeta]}\right],
    \end{align*}
    so $\frac{d\mathbb{P}^\zeta}{d\mathbb{P}}\Big|_{\widetilde{\mathcal F}_T}=\frac{\mathbbold{1}_{V_T=\zeta}}{\mathbb{P}[V_T=\zeta]}$. At $t\leq T$, for $A\in\widetilde{\mathcal F}_t$, by the tower property, we have $\frac{d\mathbb{P}^\zeta}{d\mathbb{P}}\Big|_{\widetilde{\mathcal F}_t}=\mathbb{E}^\mathbb{P}\left[\frac{d\mathbb{P}^\zeta}{d\mathbb{P}}\Big|_{\widetilde{\mathcal F}_T}\mid \widetilde{\mathcal F}_t\right]$. For ease of notation, we will use $D_t^\zeta$ to denote $\frac{d\mathbb{P}^\zeta}{d\mathbb{P}}\Big|_{\widetilde{\mathcal F}_t}$. As $V_t$ is Markov, $H_t^\zeta(V_t)=\mathbb{P}[V_T=\zeta\mid \widetilde{\mathcal{F}}_t]=\mathbb{E}^\mathbb{P}[\mathbbold{1}_{V_T=\zeta}\mid \widetilde{\mathcal F}_t]$, so $D_t^\zeta=\frac{H_t^\zeta(V_t)}{\mathbb{P}[V_T=\zeta]}$. By Bayes' Rule and the equalities $\mathbb{P}[U_0=\zeta]=f(\zeta)\pi^{\otimes N}(\zeta)$, $\mathbb{P}[U_{T-t}=x]=p_{T-t}(x)\pi^{\otimes N}(x)$, and $\mathbb{P}[U_{T-t}=x\mid U_0=\zeta]=\mathbb{E}[\mathbbold{1}_{\{x\}}(U_{T-t})\mid U_0=\zeta]=P_{T-t}^\pi\mathbbold{1}_x(\zeta)=\int_G\mathbbold{1}_x(y)\ d\pi_{T-t,\zeta}(y)=\pi_{T-t,\zeta}(x)$, we can simplify $D_t^\zeta$ to
    \begin{align*}
        D_t^\zeta=\frac{\pi_{T-t,\zeta}(V_t)}{p_{T-t}(V_t)\pi^{\otimes N}(V_t)}.
    \end{align*}
    At a specific $(t,z)$ update of $\widetilde N^{i,g}$, $V_{t^-}$ is updated to $V_{t^-}^{i\to g}$ if $0<z\leq a_{t,i,g}(V_{t^-})$ and is unchanged otherwise. Therefore
    \begin{align*}
        \frac{D_t^\zeta}{D_{t^-}^\zeta}=1+\bigl(r_{t,i,g}^\zeta(V_{t^-})-1\bigr)\mathbbold{1}_{\{0<z\leq a_{t,i,g}(V_{t^-})\}}.
    \end{align*}
    The $\mathbb P$-compensator of $\widetilde N^{i,g}$ is $\nu^{i,g}(dt,dz)=dt\,dz$. Therefore, by Jacod--Shiryaev~\cite[Theorem III.3.17]{JacodShiryaev2003}, its $\mathbb P^\zeta$-compensator is
    \begin{align*}
        \nu^{\zeta,i,g}(dt,dz)&=\left[1+\left(r_{t,i,g}^\zeta(V_{t^-})-1\right)\mathbbold{1}_{\left\{0<z\leq a_{t,i,g}(V_{t^-})\right\}}\right]dt\,dz.
    \end{align*}
    For $z\in\mathbb R$, let
    \begin{align*}
        \Gamma_{t,i,g}^{\varphi}(x,y,z)&:=\Delta_{i,g}^{xy}\varphi(x,y)\mathbbold{1}_{\{0<z\leq a_{t,i,g}(x)\wedge b_{t,i,g}(x)\}}+\Delta_{i,g}^{x}\varphi(x,y)\mathbbold{1}_{\{b_{t,i,g}(x)<z\leq a_{t,i,g}(x)\}}\\&\quad+\Delta_{i,g}^{y}\varphi(x,y)\mathbbold{1}_{\{a_{t,i,g}(x)<z\leq b_{t,i,g}(x)\}}.
    \end{align*}
    The jump equations for $(V_t,W_t)$ give that (by using the definition of the stochastic integral)
    \begin{align*}
        \varphi(V_t,W_t)-\varphi(V_\theta,W_\theta)=\sum_{i=1}^N\sum_{g\in G_0}\int\limits_\theta^t\int_{\mathbb R}\Gamma_{s,i,g}^{\varphi}(V_{s^-},W_{s^-},z)\widetilde N^{i,g}(ds,dz).
    \end{align*}
    Under $\mathbb P^\zeta$, decompose
    \begin{align*}
        \widetilde N^{i,g}(ds,dz)=\left(\widetilde N^{i,g}(ds,dz)-\nu^{\zeta,i,g}(ds,dz)\right)+\nu^{\zeta,i,g}(ds,dz).
    \end{align*}
    By the definition of a compensator,
    \begin{align*}
        \sum_{i=1}^N\sum_{g\in G_0}\int\limits_\theta^t\int_{\mathbb R}\Gamma_{s,i,g}^{\varphi}(V_{s^-},W_{s^-},z)\left(\widetilde N^{i,g}(ds,dz)-\nu^{\zeta,i,g}(ds,dz)\right)
    \end{align*}
    is a $\mathbb P^\zeta$-martingale. The integrands and jump rates are bounded on $[\theta,T_0]$ since $T-T_0=\tau>0$ and the state space is finite. Consequently, the predictable drift is
    \begin{align*}
        \mathcal G_t^\zeta\varphi(x,y)=\sum_{i=1}^N\sum_{g\in G_0}\int_{\mathbb R}\Gamma_{t,i,g}^{\varphi}(x,y,z)Y_{t,i,g}^\zeta(z)\,dz,
    \end{align*}
    where (by following Jacod), $Y_{t,i,g}^\zeta(z)=\mathbb{E}^\mathbb{P}_{\widetilde{N}^{i,g}}\left[\frac{D_t^\zeta}{D_{t^-}^\zeta}\right](t,z)$. On the region where both $V_t$ and $W_t$ are updated or where only $V_t$ is updated, the condition $z\leq a_{t,i,g}(x)$ holds, and hence $Y_{t,i,g}^\zeta(z)=r_{t,i,g}^\zeta(x)$. On the other hand, on the region where only $W_t$ is updated, $z>a_{t,i,g}(x)$, and hence $Y_{t,i,g}^\zeta(z)=1$. Integrating with respect to $z$ then gives the desired $\mathcal{G}_t^\zeta$.
\end{proof}

\begin{lemma}\label{lem:localized-comparison}
    Suppose that $V_\theta=W_\theta$, and let $A\in\widetilde{\mathcal F}_\theta$. Then, for every $\varphi:G\to[0,1]$,
    \begin{align*}
        &\mathbb E\left[\mathbbold 1_A\left(\varphi(W_{T_0})-\varphi(V_{T_0})\right)\right]=\mathbb E\left[\mathbbold 1_A\int_\theta^{T_0}\mathcal A_t^{V_T}q_t^{V_T}(V_{t^-},W_{t^-})\,dt\right],
    \end{align*}
    where $\mathcal A_t^\zeta:=\mathcal G_t^\zeta-\mathcal L_t^\zeta$.
\end{lemma}
\begin{proof}
    Fix $\zeta$ and apply the jump It\^o formula to $q_t^\zeta(V_t,W_t)$ under $\mathbb P^\zeta$. Lemma~\ref{lem:conditioned-joint-generator} gives
    \begin{align*}
        q_{T_0}^\zeta(V_{T_0},W_{T_0})-q_\theta^\zeta(V_\theta,W_\theta)&=\int_\theta^{T_0}(\partial_t+\mathcal G_t^\zeta)q_t^\zeta(V_{t^-},W_{t^-})\,dt+M_{T_0},
    \end{align*}
    where $M_t$ is a $\mathbb{P}^\zeta$-martingale, $M_\theta=0$. Multiply by $\mathbbold{1}_A$ and take conditional expectations under $\mathbb{P}^\zeta$ given $\widetilde{\mathcal F}_\theta$. Since $(\partial_t+\mathcal{L}_t^\zeta)q_t^\zeta=0$ and $\mathcal{A}_t^\zeta=\mathcal{G}_t^\zeta-\mathcal{L}_t^\zeta$, $(\partial_t+\mathcal{G}_t^\zeta)q_t^\zeta=\mathcal{A}_t^\zeta q_t^\zeta$. We also notice that $q_{T_0}^\zeta(V_{T_0},W_{T_0})=\varphi(W_{T_0})$. Since $V_\theta=W_\theta$, Lemma~\ref{lem:synchronous-conditioned} and the conditional Markov property of $V$ in the joint filtration under $\mathbb P^\zeta$ give $q_\theta^\zeta(V_\theta,W_\theta)=\mathbb E^\zeta[\varphi(V_{T_0})\mid\widetilde{\mathcal F}_\theta]$.
    Using $A\in\widetilde{\mathcal F}_\theta$ and the tower property then gives
    \begin{align*}
        \mathbb{E}[\mathbbold{1}_A(\varphi(W_{T_0})-\varphi(V_{T_0}))\mid \widetilde{\mathcal F}_\theta,V_T=\zeta]=\mathbb{E}\left[\mathbbold{1}_A\int\limits_\theta^{T_0}\mathcal{A}_t^\zeta q_t^\zeta(V_{t^-},W_{t^-})\ dt\Big|\ \widetilde{\mathcal F}_\theta,V_T=\zeta\right].
    \end{align*}
    To remove the conditions, we average over $V_T$ and then over $\widetilde{\mathcal F}_\theta$, using the tower property. This gives the desired equality.
\end{proof}

In the Boolean setting, Lu-Guo-Fang~\cite[Lemma 4.1, Equations (25)--(27)]{lu2026weaktype} use the fact that a probability measure on $\sfrac{\mathbb{Z}}{2\mathbb{Z}}$ is determined by its mean. Since this is no longer true in a state space with more than two elements, we instead have to find and then compare the full conditional coordinate laws. The following estimate controls both of the conditional law differences we will need for our proof. This definition of $\mu_t^{a,b,z}$ is motivated in the proof of Lemma~\ref{lem:coupling-energies}.

\begin{lemma}\label{lem:one-coordinate-bridge}
    Fix $z\in G_0$ and define
    \begin{align*}
        h_s(v):=e^{-(T-s)}\mathbf{1}_{\{v=z\}}+(1-e^{-(T-s)})\pi(z),\qquad 0\leq s\leq T_0,\quad v\in G_0.
    \end{align*}
    For $0\leq t<T_0$, set $r:=e^{-(T_0-t)}$ and, for $a,b,u\in G_0$, define
    \begin{align*}
        \mu_t^{a,b,z}(u):=\frac{rh_{T_0}(a)\mathbf{1}_{\{u=b\}}+(1-r)\pi(u)h_{T_0}(u)}{h_t(a)}.
    \end{align*}
    Then, for every $a,b\in G_0$,
    \begin{align*}
        \sum_{g\in G_0}\pi(g)\frac{h_t(g)}{h_t(a)}\sum_{u\in G_0}\frac{\bigl(\mu_t^{a,b,z}(u)-\mu_t^{g,b,z}(u)\bigr)^2+\bigl(\mu_t^{a,g,z}(u)-\mu_t^{g,g,z}(u)\bigr)^2}{\mu_t^{a,b,z}(u)}\leq |G_0|e^{-(T_0-t)}\kappa_{\tau,\pi}.
    \end{align*}
\end{lemma}
\begin{proof}
    For $v\in G_0$, using the equivalence $rh_{T_0}(a)+(1-r)\pi(z)=re^{-(T-T_0)}\mathbbold{1}_{a=z}+(r(1-e^{-(T-T_0)})+1-r)\pi(z)=e^{-(T-t)}\mathbbold{1}_{a=z}+(1-e^{-(T-t)})\pi(z)=h_t(a)$, we find
    \begin{align*}
        &\mu_{t}^{a,v,z}(u)-\mu_{t}^{g,v,z}(u)=\left(\frac{1}{h_t(a)}-\frac{1}{h_t(g)}\right)(1-r)(\pi(u)h_{T_0}(u)-\pi(z)\mathbbold{1}_{u=v}).
    \end{align*}
    Now consider $h_t(g)-h_t(a)=e^{-(T-t)}\mathbbold{1}_{g=z}+(1-e^{-(T-t)})\pi(z)-(e^{-(T-t)}\mathbbold{1}_{a=z}+(1-e^{-(T-t)})\pi(z))$. This is zero unless $a\neq g=z$ or $g\neq a=z$. We first consider the case where $a\neq g=z$. Here, the sum reduces to
    \begin{align*}
        &\pi(z)\frac{(h_t(z)-h_t(a))^2}{h_t(a)^2h_t(g)}(1-r)^2\sum\limits_{u\in G_0}\frac{\left(\pi(u)h_{T_0}(u)-\pi(z)\mathbbold{1}_{u=z}\right)^2+\left(\pi(u)h_{T_0}(u)-\pi(z)\mathbbold{1}_{u=b}\right)^2}{rh_{T_0}(a)\mathbbold{1}_{u=b}+(1-r)\pi(u)h_{T_0}(u)}.
    \end{align*}
    The inner sum will split into $u=b$ and $u\neq b$ terms. Noting that $rh_{T_0}(a)+(1-r)\pi(b)h_{T_0}(b)$ is greater than or equal to each individual part of it, we can bound the $b$ term by either piece of the sum. Also $|\pi(b)h_{T_0}(b)-\pi(z)|=|\pi(b)(e^{-(T-T_0)}\mathbbold{1}_{b=z}+(1-e^{-(T-T_0)})\pi(z))-\pi(z)|=\pi(z)|e^{-(T-T_0)}\mathbbold{1}_{b=z}+\pi(b)(1-e^{-(T-T_0)})-1|$, so noting that $0\leq e^{-(T-T_0)}\mathbbold{1}_{b=z}+\pi(b)(1-e^{-(T-T_0)})\leq 1$, we can bound the absolute value from above by $1$, so this is $\leq \pi(z)$. This lets us bound the sum by (recombining the matching $u=b$ term $\frac{(\pi(b)h_{T_0}(b)-\pi(z)\mathbbold{1}_{b=z})^2}{(1-r)\pi(b)h_{T_0}(b)}$ with the sum over $G_0\backslash\{b\}$ of the same quantity, then using $\sum\limits_{u\in G_0}\pi(u)h_{T_0}(u)=\pi(z)$ and expanding the square)
    \begin{align*}
        &\leq\frac{(\pi(z))^2}{rh_{T_0}(a)}+\frac{\pi(z)}{(1-r)h_{T_0}(z)}.
    \end{align*}
    Therefore the sum over $G_0$ we are bounding is (in this case) bounded by (using $e^{-(T-t)}=re^{-\tau}$ and expanding the $h_t$ terms)
    \begin{align*}
        &\leq\frac{1}{(1-re^{-\tau})^2(re^{-\tau}+(1-re^{-\tau})\pi(z))}\left(\frac{re^{-2\tau}(1-r)^2}{(1-e^{-\tau})}+\frac{r^2e^{-2\tau}(1-r)}{(e^{-\tau}+(1-e^{-\tau})\pi(z))}\right).
    \end{align*}
    We now note that $\frac{(1-r)^2}{(1-re^{-\tau})^2}\leq 1$ and $\frac{r(1-r)}{(1-re^{-\tau})^2}=\frac{r(1-r)}{(1-r+r(1-e^{-\tau}))^2}\leq\frac{r(1-r)}{4(1-r)r(1-e^{-\tau})}=\frac{1}{4(1-e^{-\tau})}$. Also, $re^{-\tau}+(1-re^{-\tau})\pi(z)=\pi(z)+re^{-\tau}(1-\pi(z))\geq\pi(z)$ and $e^{-\tau}+(1-e^{-\tau})\pi(z)\geq e^{-\tau}$. This lets us reduce the bound to simply
    \begin{align*}
        &\leq\frac{re^{-2\tau}}{(1-e^{-\tau})\pi(z)}\left(1+\frac{1}{4 e^{-\tau}}\right)=\frac{r}{(1-e^{-\tau})\pi(z)}\left(e^{-2\tau}+\frac{e^{-\tau}}{4}\right)\leq\frac{5}{4}\frac{re^{-\tau}}{(1-e^{-\tau})\pi(z)}.
    \end{align*}
    Now since $\frac{e^{-\tau}}{(1-e^{-\tau})\pi(z)}\leq\frac{e^{-\tau}}{(1-e^{-\tau})\min\limits_{g\in G_0}\pi(g)}\leq\frac{e^{-\tau}+(1-e^{-\tau})\max\limits_{g\in G_0}\pi(g)}{(1-e^{-\tau})\min\limits_{g\in G_0}\pi(g)}=\kappa_{\tau,\pi}$, we can bound the sum by $\frac{5}{4}e^{-(T_0-t)}\kappa_{\tau,\pi}$.

    We now consider the case where $g\neq a=z$. Here the sum we are interested in is
    \begin{align*}
        \sum\limits_{g\in G_0\backslash\{a\}}\frac{\pi(g)e^{-2(T-t)}}{h_t(a)^2h_t(g)}(1-r)^2\sum\limits_{u\in G_0}\frac{(\pi(u)h_{T_0}(u)-\pi(a)\mathbbold{1}_{u=g})^2+(\pi(u)h_{T_0}(u)-\pi(a)\mathbbold{1}_{u=b})^2}{rh_{T_0}(a)\mathbbold{1}_{u=b}+(1-r)\pi(u)h_{T_0}(u)}.
    \end{align*}
    We then do the same splitting, bounding, and recombining of the $u$ sum as we did in the other case. This gives a bound for the overall sum of
    \begin{align*}
        &\leq\sum\limits_{g\in G_0\backslash\{a\}}\left(\frac{re^{-2\tau}(1-r)}{(re^{-\tau}+(1-re^{-\tau})\pi(a))^2(1-re^{-\tau})}\right)\Bigg(\frac{\pi(g)(1-r)\pi(a)}{(e^{-\tau}+(1-e^{-\tau})\pi(a))}+\frac{r}{(1-e^{-\tau})}\Bigg).
    \end{align*}
    We now note that $(re^{-\tau}+(1-re^{-\tau})\pi(a))^2\geq re^{-\tau}\pi(a)$ or $\geq(\pi(a))^2$. Using this separately on both of the fractions in the second quantity and noting $\frac{(1-r)^2}{1-re^{-\tau}}\leq1-r\leq1$ and $\frac{1-r}{1-re^{-\tau}}\leq1$, we find a bound of
    \begin{align*}
        &\leq\frac{re^{-2\tau}}{\pi(a)(e^{-\tau}+(1-e^{-\tau})\pi(a))}+\frac{re^{-\tau}(|G_0|-1)}{(1-e^{-\tau})\pi(a)}\\
    \end{align*}
    Now $e^{-\tau}+(1-e^{-\tau})\pi(a)\geq e^{-\tau}$, so we have the bound
    \begin{align*}
        &\leq\frac{re^{-2\tau}}{\pi(a)e^{-\tau}}+\frac{re^{-\tau}(|G_0|-1)}{(1-e^{-\tau})\pi(a)}\leq \frac{re^{-\tau}\left(|G_0|-e^{-\tau}\right)}{\pi(a)(1-e^{-\tau})}\leq (|G_0|-e^{-\tau})r\kappa_{\tau,\pi}\leq |G_0|e^{-(T_0-t)}\kappa_{\tau,\pi}.
    \end{align*}
    Therefore in either case, the sum is bounded by $|G_0|e^{-(T_0-t)}\kappa_{\tau,\pi}$.
\end{proof}

\begin{lemma}\label{lem:coupling-energies}
    Let $\varphi:G\to[0,1]$, and let
    $A\in\widetilde{\mathcal F}_\theta$ satisfy
    $A\subseteq\mathcal E_\theta$. Define
    \begin{align*}
        \Psi_{\mathrm{sync},A}&:=\mathbb E\left[\mathbbold 1_A\int\limits_\theta^\sigma\sum_{i=1}^N\sum_{g\in G_0}\pi(g)\lambda_{t,i,g}^{V_T}(V_{t^-})\left|\Delta_{i,g}^{xy}q_t^{V_T}(V_{t^-},W_{t^-})\right|^2\,dt\right],\\
        \Psi_{\mathrm{hor},A}&:=\mathbb E\Bigg[\mathbbold 1_A\int\limits_\theta^\sigma\sum_{i=1}^N\sum_{g\in G_0}\pi(g)\lambda_{t,i,g}^{V_T}(V_{t^-})\left(\left|D(V_{t^-},W_{t^-})\right|^2+\left|D(V_{t^-},W_{t^-}^{i\to g})\right|^2\right)\,dt\Bigg],
    \end{align*}
    using the abbreviation $D$ for $\Delta_{i,g}^xq_t^{V_T}$. Then
    \begin{align*}
        \Psi_{\mathrm{hor},A}&\lesssim_\pi\kappa_{\tau,\pi}\mathbb P(A),\\
        \Psi_{\mathrm{sync},A}+\Psi_{\mathrm{hor},A}&\lesssim_\pi\kappa_{\tau,\pi}\left(\mathbb P(A)+\Lambda_{\tau,\pi}^{2}\mathscr S_A\right).
    \end{align*}
\end{lemma}
\begin{proof}
    The proof has five main steps. For the $\Psi_{\mathrm{hor},A}$ bound, we first derive a formula for the conditional coordinate law $\widehat{\mathbb{P}}^\zeta_{t,x,y}[W_{T_0}^{[i]}=u]$, then we apply Lemma~\ref{lem:one-coordinate-bridge}, then finally use Lemma~\ref{lem:rep} to bound the sum over dimensions without a factor of $N$. We then use Dynkin's formula to get a bound for $\Psi_{\mathrm{sync},A}$ in terms of $\mathcal{A}_t^\zeta(q_t^\zeta)^2$, then finally use the AM-GM inequality to bound the generator term.
    
    We begin by noticing that
    \begin{align*}
        \lambda_{t,i,g}^\zeta(x)=\frac{h_{t,i}^\zeta(g)}{h_{t,i}^\zeta(x_i)},\qquad h_{t,i}^\zeta(a):=e^{-(T-t)}\mathbbold{1}_{\{a=\zeta_i\}}+(1-e^{-(T-t)})\pi(\zeta_i).
    \end{align*}
    This follows by using the formula for $H_t^\zeta$ and canceling the density ratio (the term coming from expanding $\rho$ and one of the parts of $H_t^\zeta$). Under $\widehat{\mathbb P}_{t,x,y}^\zeta$, the coordinates of $V_t$ and $W_t$ are independent, since their rates depend only on their own coordinates. Define
    \begin{align*}
        \mu_{t,i}^{x,y,\zeta}(u):=\widehat{\mathbb P}_{t,x,y}^\zeta[\widehat W_{T_0}^{[i]}=u],\qquad\mu_t^{x,y,\zeta}:=\bigotimes_{i=1}^N\mu_{t,i}^{x,y,\zeta}.
    \end{align*}
    The coordinate law depends on $x,y,\zeta$ only through $x_i,y_i,\zeta_i$. We use shorthand that $\mu_{t,i}^{g,y,\zeta}=\mu_{t,i}^{x^{i\to g},y,\zeta}$, and $\mu_{t,i}^{a,b,z}$ denotes the law with coordinate values $a,b,z$. We also use $\mu(u)=\mu(\{u\})$. For the fixed baseline law $\mu_t^{x,y,\zeta}$, set $\bar{\varphi}:=\mathbb E^{\mu_t^{x,y,\zeta}}[\varphi]$ and
    \begin{align*}
        \varphi_i(u)&:=\sum_{(a_j)_{j\ne i}\in G_0^{N-1}}\varphi(a_1,\ldots,a_{i-1},u,a_{i+1},\ldots,a_N)\prod_{j\ne i}\mu_{t,j}^{x,y,\zeta}(a_j).
    \end{align*}
    $\varphi_i$ depends on $t,x,y,\zeta$. Since replacing $x_i$ does not change the other coordinates,
    \begin{align*}
        \Delta_{i,g}^xq_t^\zeta(x,y)&=\sum_{u\in G_0}\bigl(\mu_{t,i}^{g,y,\zeta}(u)-\mu_{t,i}^{x,y,\zeta}(u)\bigr)(\varphi_i(u)-\bar\varphi),
    \end{align*}
    where we can add in the $\bar{\varphi}$ term as the sum of the difference in laws times any constant is zero. By Cauchy-Schwarz (noting that the conditional coordinate laws have full support for $t<T_0$), we can bound $|\Delta_{i,g}^xq_t^\zeta(x,y)|^2$ and $|\Delta_{i,g}^xq_t^\zeta(x,y^{i\to g})|^2$ from above by the following respectively:
    \begin{align*}
        &\left(\sum\limits_{u\in G_0}\left(\frac{\mu_{t,i}^{x,y,\zeta}(u)-\mu_{t,i}^{g,y,\zeta}(u)}{\sqrt{\mu_{t,i}^{x,y,\zeta}(u)}}\right)^2\right)\left(\sum\limits_{u\in G_0}\left(\sqrt{\mu_{t,i}^{x,y,\zeta}(u)}(\varphi_i(u)-\bar{\varphi})\right)^2\right),\\
        &\left(\sum\limits_{u\in G_0}\left(\frac{\mu_{t,i}^{x,g,\zeta}(u)-\mu_{t,i}^{g,g,\zeta}(u)}{\sqrt{\mu_{t,i}^{x,y,\zeta}(u)}}\right)^2\right)\left(\sum\limits_{u\in G_0}\left(\sqrt{\mu_{t,i}^{x,y,\zeta}(u)}(\varphi_i(u)-\bar{\varphi})\right)^2\right).
    \end{align*}
    This allows us to bound $\Psi_{\mathrm{hor},A}$ by
    \begin{align*}
        \Psi_{\mathrm{hor},A}\leq\mathbb{E}\Bigg[\mathbbold{1}_A\int\limits_\theta^\sigma\sum\limits_{i=1}^N\left(\sum\limits_{u\in G_0}\left(\sqrt{\mu_{t,i}^{V_{t^-},W_{t^-},{V_T}}(u)}(\varphi_i(u)-\bar{\varphi})\right)^2\right)\sum\limits_{g\in G_0}\pi(g)\lambda_{t,i,g}^{V_T}(V_{t^-})\times\\\left(\sum\limits_{u\in G_0}\left(\frac{\mu_{t,i}^{V_{t^-},g,{V_T}}(u)-\mu_{t,i}^{g,g,{V_T}}(u)}{\sqrt{\mu_{t,i}^{V_{t^-},W_{t^-},{V_T}}(u)}}\right)^2+\sum\limits_{u\in G_0}\left(\frac{\mu_{t,i}^{V_{t^-},W_{t^-},{V_T}}(u)-\mu_{t,i}^{g,W_{t^-},{V_T}}(u)}{\sqrt{\mu_{t,i}^{V_{t^-},W_{t^-},{V_T}}(u)}}\right)^2\right)\ dt\Bigg].
    \end{align*}
    We now consider the inner sum over $g\in G_0$. To do this, we first expand $\mu_{t,i}^{x,y,\zeta}(u)$. We note that as in Lemma~\ref{lem:synchronous-conditioned}, $(\widehat{V}_{t}^{[i]},\widehat{W}_{t}^{[i]})$ has joint generator
    \begin{align*}
        \mathcal{L}_t^{\zeta_i}f(a,b)=\sum\limits_{w\in G_0}\pi(w)(f(w,w)-f(a,b))\frac{h_t(w)}{h_t(a)},
    \end{align*}
    where we fix $i,\zeta$, and set $h_t(x)=e^{-(T-t)}\mathbbold{1}_{\{x=\zeta_i\}}+(1-e^{-(T-t)})\pi(\zeta_i)$ as in Lemma~\ref{lem:one-coordinate-bridge}. This can be seen by simply computing the probabilities in the one-dimensional version of $\lambda_{t,i,g}^\zeta$ in the Lemma~\ref{lem:synchronous-conditioned} statement. By definition and Lemma~\ref{lem:synchronous-conditioned}, $\mu_{t,i}^{x,y,\zeta}$ will satisfy the PDE $(\partial_t+\mathcal{L}_t^{\zeta_i})\mu_{t,i}^{x,y,\zeta}(u)=0$. Multiplying both sides by $h_t(x_i)$ and expanding the generator then gives
    \begin{align*}
        -h_t(x_i)\partial_t\mu_{t,i}^{x,y,\zeta}(u)        &=\sum\limits_{w\in G_0}\pi(w)\mu_{t,i}^{w,w,\zeta}(u)h_t(w)-\mu_{t,i}^{x,y,\zeta}(u)\sum\limits_{w\in G_0}\pi(w)h_t(w).
    \end{align*}
    By definition of $h_t$, the final sum is
    \begin{align*}
        \sum_{w\in G_0}\pi(w)h_t(w)=e^{-(T-t)}\pi(\zeta_i)+(1-e^{-(T-t)})\pi(\zeta_i)=\pi(\zeta_i).
    \end{align*}
    
    The first one can have the $\mu h$ product rewritten as $M_t(w,w,u)$ where $M_t(x,y,u):=h_t(x_i)\mu_{t,i}^{x,y,\zeta}(u)$. We use shorthand of $M_t(a,b,u):=h_t(a)\mu_{t,i}^{a,b,\zeta}(u)$ for $a,b\in G_0$. Thus
    \begin{align*}
        -h_t(x_i)\partial_t\mu_{t,i}^{x,y,\zeta}(u)=\sum\limits_{w\in G_0}\pi(w)M_t(w,w,u)-\pi(\zeta_i)\mu_{t,i}^{x,y,\zeta}(u).
    \end{align*}
    Now consider $\partial_tM_t=\partial_th_t\mu_{t,i}+\partial_t\mu_{t,i}h_t$. This simplifies to
    \begin{align*}
        \partial_tM_t(x,y,u)&=M_t(x,y,u)-\sum\limits_{w\in G_0}\pi(w)M_t(w,w,u).
    \end{align*}
    We now consider $\frac{d}{dt}\left(\sum\limits_{w\in G_0}\pi(w)M_t(w,w,u)\right)$. By the equation above, this simplifies to (as $\pi$ is a probability measure)
    \begin{align*}
        \sum\limits_{w\in G_0}\pi(w)M_t(w,w,u)-\sum\limits_{w\in G_0}\pi(w)\sum\limits_{w'\in G_0}\pi(w')M_t(w',w',u)=0.
    \end{align*}
    This means that this term is constant with respect to time, so we can just set $t=T_0$ in that term. We evaluate 
    \begin{align*}
        &\sum\limits_{w\in G_0}\pi(w)M_{T_0}(w,w,u)=\sum\limits_{w\in G_0}\pi(w)(e^{-(T-T_0)}\delta_{w=\zeta_i}+(1-e^{-(T-T_0)})\pi(\zeta_i))\mathbbold{1}_{u=w}\\
        &=\pi(u)(e^{-(T-T_0)}\delta_{u=\zeta_i}+(1-e^{-(T-T_0)})\pi(\zeta_i))=\pi(u)h_{T_0}(u).
    \end{align*}
    This gives $\partial_tM_t(x,y,u)=M_t(x,y,u)-\pi(u)h_{T_0}(u)$. This is an ODE which we can solve as we know the terminal condition $M_{T_0}$. The solution is
    \begin{align*}
        M_t(x,y,u)=e^{-(T_0-t)}h_{T_0}(x_i)\mathbbold{1}_{u=y_i}+(1-e^{-(T_0-t)})\pi(u)h_{T_0}(u).
    \end{align*}
    This then lets us simplify
    \begin{equation}\label{eq:coordinate-bridge}
        \mu_{t,i}^{x,y,\zeta}(u)=\frac{e^{-(T_0-t)}h_{T_0}(x_i)\mathbbold{1}_{u=y_i}+(1-e^{-(T_0-t)})\pi(u)h_{T_0}(u)}{h_t(x_i)}.
    \end{equation}
    Now, by Lemma~\ref{lem:one-coordinate-bridge}, applied with $a=(V_{t^-})_i,b=(W_{t^-})_i,z=(V_T)_i,\mu_t^{a,b,z}=\mu_{t,i}^{a,b,z}$, we have that the inner sum over $G_0$ is bounded by $|G_0|e^{-(T_0-t)}\kappa_{\tau,\pi}$. This gives the $\Psi_{\mathrm{hor},A}$ bound
    \begin{align*}
        \Psi_{\mathrm{hor},A}\leq\mathbb{E}\Bigg[\mathbbold{1}_A\int\limits_\theta^\sigma\sum\limits_{i=1}^N\left(\sum\limits_{u\in G_0}\left(\sqrt{\mu_{t,i}^{V_{t^-},W_{t^-},{V_T}}(u)}(\varphi_i(u)-\bar{\varphi})\right)^2\right)|G_0|e^{-(T_0-t)}\kappa_{\tau,\pi}\ dt\Bigg].
    \end{align*}
    The double sum is bounded by $\frac14$ by Lemma~\ref{lem:rep}. Thus
    \begin{align*}
        \Psi_{\mathrm{hor},A}\leq\mathbb{E}\Bigg[\mathbbold{1}_A\int\limits_\theta^\sigma\frac{|G_0|}{4}e^{-(T_0-t)}\kappa_{\tau,\pi}\ dt\Bigg]\leq\frac{|G_0|\kappa_{\tau,\pi}}{4}\mathbb{E}[\mathbbold{1}_A(\sigma-\theta)]\leq\frac{|G_0|\kappa_{\tau,\pi}}{4}\mathbb{P}[A].
    \end{align*}
    This gives our desired bound, where the $\pi$ dependence in the constant is through $|G_0|=|\mathrm{supp}(\pi)|$.

    We now consider $\Psi_{\mathrm{sync},A}$. For each fixed $\zeta$, apply time-dependent Dynkin's formula under $\mathbb P^\zeta$ to $(q_t^\zeta)^2$, with process $(V_t,W_t)$, stopping time $\sigma$, and initialization time $\theta$, and then average over $V_T=\zeta$ and use the tower property to transition the expectations into $\mathbb{P}$ expectations, not $\mathbb{P}^\zeta$ ones. This gives
    \begin{align*}
        &\mathbb{E}[\mathbbold{1}_A(q_{\sigma}^{V_T}(V_{\sigma},W_{\sigma}))^2]=\mathbb{E}[\mathbbold{1}_A(q_\theta^{V_T}(V_\theta,W_\theta))^2]+\mathbb{E}\left[\mathbbold{1}_A\int\limits_\theta^{\sigma}(\partial_t+\mathcal{G}_t^{V_T})(q_t^{V_T})^2(V_{t^-},W_{t^-})\ dt\right].
    \end{align*}
    We now expand $(\partial_t+\mathcal{L}_t^\zeta)(q_t^\zeta)^2=-2q_t^\zeta\mathcal{L}_t^\zeta(q_t^\zeta)+\mathcal{L}_t^\zeta(q_t^\zeta)^2$. Using the generator formula for $\mathcal{L}_t^\zeta$ and factoring, we can simplify this to
    \begin{align*}
        &(\partial_t+\mathcal{L}_t^\zeta)(q_t^\zeta)^2=\sum\limits_{i=1}^N\sum\limits_{g\in G_0}\pi(g)\lambda_{t,i,g}^\zeta(x)|\Delta_{i,g}^{xy}q_t^\zeta(x,y)|^2.
    \end{align*}
    Therefore
    \begin{align*}
        &\Psi_{\mathrm{sync},A}\leq\mathbb{P}[A]+\mathbb{E}\left[\mathbbold{1}_A\int\limits_\theta^\sigma|\mathcal{A}_t^{V_T}(q_t^{V_T})^2(V_{t^-},W_{t^-})|\ dt\right].
    \end{align*}
    We now notice that
    \begin{align*}
        \Delta_{i,g}^y\varphi(x^{i\to g},y)=\Delta_{i,g}^{xy}\varphi(x,y)-\Delta_{i,g}^x\varphi(x,y),\qquad\Delta_{i,g}^y\varphi(x,y)   =\Delta_{i,g}^{xy}\varphi(x,y)-\Delta_{i,g}^x\varphi(x,y^{i\to g}).
    \end{align*}
    This lets us bound $|\mathcal{A}_t^\zeta(q_t^\zeta)^2|$ by (using $|\Delta_{i,g}^y(q_t^\zeta)^2|\leq 2|\Delta_{i,g}^yq_t^\zeta|$): $|\mathcal{A}_t^\zeta(q_t^\zeta)^2(x,y)|\leq$
    \begin{align*}        &2\sum\limits_{i=1}^N\sum\limits_{g\in G_0}\Big[(r_{t,i,g}^\zeta(x)(a_{t,i,g}(x)-b_{t,i,g}(x))_++(b_{t,i,g}(x)-a_{t,i,g}(x))_+)|\Delta_{i,g}^{xy}q_t^\zeta(x,y)|\\&+r_{t,i,g}^\zeta(x)(a_{t,i,g}(x)-b_{t,i,g}(x))_+|\Delta_{i,g}^xq_t^\zeta(x,y)|+(b_{t,i,g}(x)-a_{t,i,g}(x))_+|\Delta_{i,g}^xq_t^\zeta(x,y^{i\to g})|\Big].
    \end{align*}
    We now abbreviate $r_{t,i,g}^\zeta(x)(a_{t,i,g}(x)-b_{t,i,g}(x))_+=c_+$ and $(b_{t,i,g}(x)-a_{t,i,g}(x))_+=c_-$. Inside expectations these coefficients are evaluated at $x=V_{t^-}$ and $\zeta=V_T$. We then use the AM-GM inequality separately to the two $\Delta^x$ terms to get 
    \begin{align*}
        2c_+|\Delta_{i,g}^xq_t^\zeta(x,y)|&\leq\pi(g)\lambda_{t,i,g}^\zeta(x)|\Delta_{i,g}^xq_t^\zeta(x,y)|^2+\frac{(c_+)^2}{\pi(g)\lambda_{t,i,g}^\zeta(x)},\\
        2c_-|\Delta_{i,g}^xq_t^\zeta(x,y^{i\to g})|&\leq\pi(g)\lambda_{t,i,g}^\zeta(x)|\Delta_{i,g}^xq_t^\zeta(x,y^{i\to g})|^2+\frac{(c_-)^2}{\pi(g)\lambda_{t,i,g}^\zeta(x)}.
    \end{align*}
    Therefore
    \begin{align*}
        &\Psi_{\mathrm{sync},A}\leq\mathbb{P}[A]+\Psi_{\mathrm{hor},A}+\mathbb{E}\Bigg[\mathbbold{1}_A\int\limits_\theta^\sigma\sum\limits_{i=1}^N\sum\limits_{g\in G_0}2(c_++c_-)|\Delta_{i,g}^{xy}q_t^{V_T}(V_{t^-},W_{t^-})|\ dt+\\&\mathbbold{1}_A\int\limits_\theta^\sigma\sum\limits_{i=1}^N\sum\limits_{g\in G_0}\Big[\frac{(c_+)^2+(c_-)^2}{\pi(g)\lambda_{t,i,g}^{V_T}(V_{t^-})}\Big]\ dt\Bigg].
    \end{align*}
    We now apply AM-GM again to the remaining term: $2(c_++c_-)|\Delta_{i,g}^{xy}q_t^\zeta(x,y)|\leq\\\frac{1}{2}\pi(g)\lambda_{t,i,g}^\zeta(x)|\Delta_{i,g}^{xy}q_t^\zeta(x,y)|^2+2\frac{(c_++c_-)^2}{\pi(g)\lambda_{t,i,g}^\zeta(x)}$. $(c_++c_-)^2=c_+^2+c_-^2$ as $c_+$ and $c_-$ are never both nonzero. Then recognizing the definition of $\Psi_{\mathrm{sync},A}$ in the AM-GM term above,
    \begin{align*}
        &\Psi_{\mathrm{sync},A}\leq\mathbb{P}[A]+\Psi_{\mathrm{hor},A}+\frac{1}{2}\Psi_{\mathrm{sync},A}+3\mathbb{E}\Bigg[\mathbbold{1}_A\int\limits_\theta^\sigma\sum\limits_{i=1}^N\sum\limits_{g\in G_0}\frac{(c_+)^2+(c_-)^2}{\pi(g)\lambda_{t,i,g}^{V_T}(V_{t^-})}\ dt\Bigg].
    \end{align*}
    Finally, we bound the fraction. For clarity, we write $\rho=\rho_{t,i,g}(V_{t^-})$, $\lambda=\lambda_{t,i,g}^\zeta(V_{t^-})$, $\delta=\delta_t$, and $\Lambda=\Lambda_{\tau,\pi}$. Since $r_{t,i,g}^\zeta=\frac{\lambda}{\rho}$, it expands to
    \begin{align*}
        \frac{c_+^2+c_-^2}{\pi(g)\lambda}=\frac{\pi(g)\lambda}{\rho^2}   (\rho-\Psi_\delta(\rho))_+^2+\frac{\pi(g)}{\lambda}(\Psi_\delta(\rho)-\rho)_+^2.
    \end{align*}
    If $\rho\leq1$, then $\Psi_\delta(\rho)-\rho=1-\rho^\delta\geq0$, so only the second term can be nonzero. If $\rho>1$, then $\rho-\Psi_\delta(\rho)=\rho(1-\rho^{-\delta})\geq0$, so only the first term can be nonzero. Thus
    \begin{align*}
        \frac{c_+^2+c_-^2}{\pi(g)\lambda}=\begin{cases}\displaystyle\frac{\pi(g)}{\lambda}(1-\rho^\delta)^2, &\rho\leq1,\\[1mm]
          \displaystyle\pi(g)\lambda(1-\rho^{-\delta})^2,&\rho>1.
        \end{cases}
    \end{align*}
    We now see that $1-\rho^\delta=1-e^{-\delta(-\log \rho)}\leq-\delta\log\rho\leq\delta\Lambda(1-\rho)$ and $1-\rho^{-\delta}=1-e^{-\delta\log\rho}\leq\delta\log\rho\leq\delta\Lambda\frac{\rho-1}{\rho}$, so we have the bound
    \begin{align*}
        &\leq\begin{cases}
            \frac{\pi(g)}{\lambda_{t,i,g}^\zeta(V_{t^-})}\Lambda_{\tau,\pi}^2\bar{\delta}^2(1-\rho_{t,i,g}(V_{t^-}))^2\ &\rho_{t,i,g}(V_{t^-})\leq 1\\
            \frac{\lambda_{t,i,g}^\zeta(V_{t^-})\pi(g)}{\rho_{t,i,g}(V_{t^-})^2}\Lambda_{\tau,\pi}^2\bar{\delta}^2(1-\rho_{t,i,g}(V_{t^-}))^2\ &\rho_{t,i,g}(V_{t^-})> 1
        \end{cases}.
    \end{align*}
    From $\lambda_{t,i,g}^\zeta(x)=\frac{h_t(g)}{h_t(x_i)}$, both $\lambda_{t,i,g}^\zeta(x)$ and its reciprocal are at most $\kappa_{\tau,\pi}$. Therefore both cases are bounded by $\pi(g)\Lambda_{\tau,\pi}^2\kappa_{\tau,\pi}\bar\delta^2(1-\rho_{t,i,g}(V_{t^-}))^2$, so
    \begin{align*}
        \frac{1}{2}\Psi_{\mathrm{sync},A}&\leq\mathbb{P}[A]+\Psi_{\mathrm{hor},A}+3\Lambda_{\tau,\pi}^2\kappa_{\tau,\pi}\mathbb{E}\Bigg[\mathbbold{1}_A\int\limits_\theta^\sigma\sum\limits_{i=1}^N\sum\limits_{g\in G_0}\pi(g)\bar\delta^2(1-\rho_{t,i,g}(V_{t^-}))^2\ dt\Bigg].
    \end{align*}
    Now using our $\Psi_{\mathrm{hor},A}$ bound,
    \begin{align*}
        \Psi_{\mathrm{sync},A}+\Psi_{\mathrm{hor},A}&\lesssim_\pi2\mathbb{P}[A]+3\kappa_{\tau,\pi}\mathbb{P}[A]+6\Lambda_{\tau,\pi}^2\kappa_{\tau,\pi}\mathscr{S}_A\lesssim_\pi\kappa_{\tau,\pi}(\mathbb{P}[A]+\Lambda_{\tau,\pi}^2\mathscr{S}_A).
    \end{align*}
\end{proof}

We can finally prove Lemma~\ref{lem:localized-tv}.

\begin{proof}[Proof of Lemma~\ref{lem:localized-tv}]
    Start with the identity in Lemma~\ref{lem:localized-comparison}. For $x=V_{t^-}$, $y=W_{t^-}$, and $\zeta=V_T$, write
    \begin{align*}
        c_+=r_{t,i,g}^\zeta(x)(a_{t,i,g}(x)-b_{t,i,g}(x))_+,
        \quad c_-=(b_{t,i,g}(x)-a_{t,i,g}(x))_+,
        \quad w=\pi(g)\lambda_{t,i,g}^\zeta(x).
    \end{align*}
    Since $c_+c_-=0$ and $\Delta_{i,g}^yq(x,y)=\Delta_{i,g}^{xy}q(x,y)-\Delta_{i,g}^xq(x,y^{i\to g})$, $\mathcal A_t^\zeta q_t^\zeta(x,y)=$
    \begin{align*}
        \sum_{i=1}^N\sum_{g\in G_0}\Bigl[          c_+\bigl(\Delta_{i,g}^xq_t^\zeta(x,y)-\Delta_{i,g}^{xy}q_t^\zeta(x,y)\bigr)+c_-\bigl(\Delta_{i,g}^{xy}q_t^\zeta(x,y)-\Delta_{i,g}^xq_t^\zeta(x,y^{i\to g})\bigr)\Bigr].
    \end{align*}
    Both coefficients vanish after $\sigma$. Define
    \begin{align*}
        \mathcal C_A:=\mathbb E\left[\mathbbold{1}_A\int_\theta^\sigma\sum_{i=1}^N\sum_{g\in G_0}\frac{c_+^2+c_-^2}{w}\,dt\right].
    \end{align*}
    Cauchy--Schwarz on the Lemma~\ref{lem:localized-comparison} identity applied jointly over time, coordinates, branches, and probability, followed by $|u-v|^2\leq2u^2+2v^2$, gives
    \begin{align*}
        \left|\mathbb E[\mathbbold{1}_A(\varphi(W_{T_0})-\varphi(V_{T_0}))]\right|\leq\sqrt{2\mathcal C_A(\Psi_{\mathrm{sync},A}+\Psi_{\mathrm{hor},A})}.
    \end{align*}
    The coefficient estimate we found in the proof of Lemma~\ref{lem:coupling-energies} gives
    \begin{equation*}\label{eq:coefficient-energy}
        \frac{c_+^2+c_-^2}{w}
        \leq\kappa_{\tau,\pi}\Lambda_{\tau,\pi}^2
          \pi(g)\bar\delta^2(\rho_{t,i,g}(V_{t^-})-1)^2,
        \qquad t\leq\sigma.
    \end{equation*}
    Thus $\mathcal C_A\leq\kappa_{\tau,\pi}\Lambda_{\tau,\pi}^2\mathscr S_A$.
    Applying Lemma~\ref{lem:coupling-energies} now yields
    \begin{align*}
        \left|\mathbb E[\mathbbold{1}_A(\varphi(W_{T_0})-\varphi(V_{T_0}))]\right|
        &\lesssim_\pi\kappa_{\tau,\pi}\Lambda_{\tau,\pi}
          \sqrt{\mathscr S_A(\mathbb P(A)+\Lambda_{\tau,\pi}^2\mathscr S_A)}\\
        &\lesssim_\pi\kappa_{\tau,\pi}\Lambda_{\tau,\pi}
          \sqrt{\mathscr S_A\mathbb P(A)}
          +\kappa_{\tau,\pi}\Lambda_{\tau,\pi}^2\mathscr S_A.
    \end{align*}
    Finally, we take the supremum over $\varphi$ to get
    \begin{align*}
        D_A\lesssim_\pi \kappa_{\tau,\pi}\Lambda_{\tau,\pi}\sqrt{\mathscr{S}_A\mathbb{P}(A)}+\kappa_{\tau,\pi}\Lambda_{\tau,\pi}^2\mathscr{S}_A,
    \end{align*}
    noting the following equivalence. Let $\mu_A^W(B):=\mathbb{P}(A\cap\{W_{T_0}\in B\})$, $\mu_A^V(B):=\mathbb{P}(A\cap\{V_{T_0}\in B\})$ and $\nu_A:=\mu_A^W-\mu_A^V$ be measures. Then $D_A=\sup\limits_{B\subseteq G}|\nu_A(B)|$. For all $\varphi$, since $0\leq\varphi\leq1$, we use the layer cake representation
    \begin{align*}
        |\mathbb{E}[\mathbbold{1}_A(\varphi(W_{T_0})-\varphi(V_{T_0}))]|=\left|\int_G\varphi\ d\nu_A\right|\leq\int\limits_0^1|\nu_A(\{\varphi>s\})|\ ds\leq\sup\limits_{B\subseteq G}|\nu_A(B)|.
    \end{align*}
    Thus $|\mathbb{E}[\mathbbold{1}_A(\varphi(W_{T_0})-\varphi(V_{T_0}))]|\leq D_A$. For the opposite direction, simply choose $\varphi=\mathbbold{1}_B$. Therefore, $\sup\limits_{0\leq\varphi\leq1}|\mathbb{E}[\mathbbold{1}_A(\varphi(W_{T_0})-\varphi(V_{T_0}))]|= D_A$.
\end{proof}

\section{Band contraction}\label{sec6}
In this section we prove Lemma~\ref{lem:strong-monotonicity}. Our proof follows the band contraction proof in Lu, Guo, and Fang~\cite{lu2026weaktype}. To do this, we first need a helper lemma.

\begin{lemma}\label{lem:super}
    Fix $h:G\to[0,\infty)$. Define
    \begin{align*}
        \Xi_t=\begin{cases}
            \log P_{T-t}^\pi f(W_t)-(1-\bar{\delta})\log P_{T-t}^\pi f(V_t)-\bar{\delta}\log P_{T-\theta}^\pi f(V_\theta)\ &t\leq\sigma\\
            \log P_{T-t}^\pi f(W_t)-\log P_{T-t}^\pi f(V_t)+\bar{\delta}(\log P_{T-\sigma}^\pi f(V_\sigma)-\log P_{T-\theta}^\pi f(V_\theta))\ &t\geq \sigma
        \end{cases},
    \end{align*}
    let $M_t=e^{\Xi_t}$, and let $H_t(x)=\mathbb{E}[h(V_{T_0})\mid V_t=x]$. Then the process $M_tH_t(W_t)$ is a nonnegative supermartingale on $[\theta,T_0]$. Consequently, for every $A\in\widetilde{\mathcal{F}}_\theta$,
    \begin{align*}
        \mathbb{E}[\mathbbold{1}_AM_{T_0}h(W_{T_0})]\leq\mathbb{E}[\mathbbold{1}_Ah(V_{T_0})].
    \end{align*}
\end{lemma}
\begin{proof}
    We first note that if $M_tH_t(W_t)$ is a nonnegative supermartingale, we will have
    \begin{align*}
        \mathbb{E}[\mathbbold{1}_AM_{T_0} H_{T_0}(W_{T_0})]\leq \mathbb{E}[\mathbbold{1}_A M_\theta H_\theta(W_\theta)].
    \end{align*}
    Then, noting that $H_{T_0}=h$, $M_\theta=1$, $V_\theta=W_\theta$, $H_\theta(V_\theta)=\mathbb{E}[h(V_{T_0})\mid \widetilde{\mathcal{F}}_\theta]$, and $A\in\widetilde{\mathcal{F}}_\theta$, we get
    \begin{align*}
        \mathbb{E}[\mathbbold{1}_AM_{T_0} h(W_{T_0})]\leq \mathbb{E}[\mathbbold{1}_A h(V_{T_0})].
    \end{align*}
    Thus it suffices to show that $M_tH_t(W_t)$ is a supermartingale. This will hold if its predictable drift is always less than or equal to zero. We note that the two formulas for $\Xi_t$ agree at $\sigma$, so changing formulas introduces no additional jump. Let $\Phi_t(x,y)=M_tH_t(y)$, with the convention that the $V_t$ in the $\Xi_t$ formula are replaced by $x$ and the $W_t$ are replaced by $y$. By Jump It\^o, there is a martingale $Q_t$ such that
    \begin{equation}\label{eq88}
        d\Phi_t=\partial_t\Phi_t(V_{t^-},W_{t^-})\ dt+\sum\limits_{i=1}^N\sum\limits_{g\in G_0}\int_{\mathbb{R}}\Phi_t(V_t^z,W_t^z)-\Phi_t(V_{t^-},W_{t^-})\ dz\ dt+dQ_t,
    \end{equation}
    where, for each $i,g$, $V_t^z$ and $W_t^z$ are the states obtained by applying the updates using $z$ to $V_{t^-}$ and $W_{t^-}$. We first consider the case where $t\leq\sigma$. For a synchronous update, $H_t(W_t)=H_{t^-}(W_{t^-}^{i\to g})$, so we have $M_t=M_{t^-}\frac{\rho_{t,i,g}(W_{t^-})}{(\rho_{t,i,g}(V_{t^-}))^{1-\bar{\delta}}}$. If only $W_t$ is updated, then we still have $H_t(W_t)=H_{t^-}(W_{t^-}^{i\to g})$, but now $M_t=M_{t^-}\rho_{t,i,g}(W_{t^-})$. Finally, if only $V_t$ is updated, we have $H_t(W_t)=H_{t^-}(W_{t^-})$, so $M_t=M_{t^-}(\rho_{t,i,g}(V_{t^-}))^{\bar{\delta}-1}$. Using this, the integral in Equation~\ref{eq88} simplifies to (using shorthand $H=H_{t^-}(W_{t^-})$, $H_i^g=H_{t^-}(W_{t^-}^{i\to g})$, $Y=\rho_{t,i,g}(W_{t^-})$, $X=\rho_{t,i,g}(V_{t^-})$, $M=M_{t^-}$, $a_{t,i}^g=\pi(g)X$, and $b_{t,i}^g=\pi(g)\Psi_{\delta_t}(X)$)
    \begin{align*}
        &\pi(g)M\left(X^{\bar{\delta}}YH_i^g+\Psi_{\delta_t}(X)YH_i^g-\Psi_{\delta_t}(X)H-XYH_i^g\right)\ &a_{t,i}^g<b_{t,i}^g\\
        &\pi(g)M\left(\Psi_{\delta_t}(X)YX^{\bar{\delta}-1}H_i^g+X^{\bar{\delta}}H-XH-\Psi_{\delta_t}(X)X^{\bar{\delta}-1}H\right)\ &a_{t,i}^g>b_{t,i}^g.
    \end{align*}
     We now note that if $a_{t,i}^g<b_{t,i}^g$, then $\rho_{t,i,g}(V_{t^-})<\Psi_{\delta_t}(\rho_{t,i,g}(V_{t^-}))$. We note that since $\bar{\delta}\in[0,1]$, $(\rho_{t,i,g}(V_{t^-}))^{1-\bar{\delta}}\leq\rho_{t,i,g}(V_{t^-})$ if $\rho_{t,i,g}(V_{t^-})>1$, so $\rho_{t,i,g}(V_{t^-})<\Psi_{\delta_t}(\rho_{t,i,g}(V_{t^-}))$ is impossible. Therefore this case corresponds to the $0<r\leq 1$ branch of $\Psi$. The other case will correspond to the $r>1$ branch. This means that the integral expressions match and simplify to $\pi(g)M\left(YH_i^g-XH-H+X^{\bar{\delta}}H\right)$. We now consider the $\partial_t\Phi_t$ term. By the chain rule, this expands to
    \begin{align*}
        \partial_t\Phi_t(V_{t^-},W_{t^-})=M\Bigl[\partial_tH+H\bigl(\partial_t\log P_{T-t}^\pi f(W_{t^-})-(1-\bar{\delta})\partial_t\log P_{T-t}^\pi f(V_{t^-})\bigr)\Bigr].
    \end{align*}
    To compute $\partial_t\log P_{T-t}^\pi f(x)$, we note that $p_t$ satisfies $\partial_tp_t=Lp_t$, so by the chain rule, $\partial_tp_{T-t}=-Lp_{T-t}$. Therefore $\partial_t\log P_{T-t}^\pi f(x)=\frac{\partial_tp_{T-t}(x)}{p_{T-t}(x)}=-\frac{Lp_{T-t}(x)}{p_{T-t}(x)}$. Then by expanding the generator, we find that $\partial_t\log P_{T-t}^\pi f(x)=\sum\limits_{i=1}^N\sum\limits_{g\in G_0}\pi(g)(1-\rho_{t,i,g}(x))$. Therefore (recalling that our shorthand $X$ and $Y$ depend on $i,g$), we have $d(M_tH_t(W_t))=$
    \begin{align*}
        M\partial_tH\ dt+M\sum\limits_{i=1}^N\sum\limits_{g\in G_0}\pi(g)\left(\bar{\delta}H(1-X)-HY+YH_i^g+X^{\bar{\delta}}H-H 
        \right)dt+dQ_t.
    \end{align*}
    By Lemma~\ref{reverse_generator}, the generator of the time-reversed process $\widetilde{L}_tH$ in this shorthand will be $\sum\limits_{i=1}^N\sum\limits_{g\in G_0}\pi(g)Y(H_i^g-H)$. We know that $(\partial_t+\widetilde{L}_t)H=0$, so we have a simplification
    \begin{align*}
        d(M_tH_t(W_t))=M\sum\limits_{i=1}^N\sum\limits_{g\in G_0}\pi(g)H\left(\bar{\delta}-\bar{\delta}X-1+X^{\bar{\delta}} 
        \right)dt+dQ_t.
    \end{align*}
    It thus suffices to show that $\pi(g)MH\left(\bar{\delta}-\bar{\delta}X-1+X^{\bar{\delta}}\right)\leq 0$ always. Clearly $\pi(g),M,\\H\geq0$. The weighted AM-GM inequality gives $1^{1-\bar{\delta}}X^{\bar{\delta}}\leq(1-\bar{\delta})+\bar{\delta}X$, so $(1-\bar{\delta})+\bar{\delta}X-X^{\bar{\delta}}\geq(1-\bar{\delta})+\bar{\delta}X-(1-\bar{\delta})-\bar{\delta}X=0$, and we are done.

    We now consider $t>\sigma$. The coupling is now synchronous, so only synchronous updates occur, which in the same shorthand above will contribute $M\pi(g)(YH_i^g-XH)$. The $\partial_t\Phi_t$ part is now just $M\partial_tH_t(W_{t^-})+MH\sum\limits_{i=1}^N\sum\limits_{g\in G_0}\pi(g)(X-Y)$, so we have
    \begin{align*}
        d(M_tH_t(W_t))=M\partial_tH_t(W_{t^-})\ dt+M\sum\limits_{i=1}^N\sum\limits_{g\in G_0}\pi(g)Y(H_i^g-H)\ dt+dQ_t.
    \end{align*}
    The drift is just zero here. Therefore for any $t$, we have non-positive drift, so $M_tH_t(W_t)$ is a nonnegative supermartingale.
\end{proof}

Our proof of Lemma~\ref{lem:strong-monotonicity} is identical to the proof of~\cite[Lemma 3.3]{lu2026weaktype}.

\begin{proof}[Proof of Lemma~\ref{lem:strong-monotonicity}]
    For $j\in\mathbb Z_{\geq0}$, on the event
    \begin{align*}
        C_j:=\mathcal{E}_\theta\cap\{\log P_{T-T_0}^\pi f(W_{T_0})\in(\ell+j,\ell+j+1]\}\cap\{\log P_{T-T_0}^\pi f(V_{T_0})\leq\ell+1\},
    \end{align*}
    we have $\Xi_{T_0}>\alpha+j-1$. Therefore on the event $C_j$, $M_{T_0}>c_0^{-1}e^j$, so
    \begin{align*}
        \mathbb{E}[\mathbbold{1}_{C_j}]\leq c_0e^{-j}\mathbb{E}[\mathbbold{1}_{\mathcal{E}_\theta}M_{T_0}\mathbbold{1}_{\log P_{T-T_0}^\pi f(W_{T_0})\in(\ell+j,\ell+j+1]}].
    \end{align*}
    Then by Lemma~\ref{lem:super} applied on the right side with $A=\mathcal{E}_\theta$ and $\\h(x)=\mathbbold{1}_{\log P_{T-T_0}^\pi f(x)\in(\ell+j,\ell+j+1]}$,
    \begin{align*}
        \mathbb{E}[\mathbbold{1}_{C_j}]&\leq c_0e^{-j}E[\mathbbold{1}_{\mathcal{E}_\theta}\mathbbold{1}_{\log P_{T-T_0}^\pi f(V_{T_0})\in(\ell+j,\ell+j+1]}]\\
        \mathbb{P}[C_j]&\leq c_0e^{-j}\mathbb{P}[\mathcal{E}_\theta\cap\{\log P_{T-T_0}^\pi f(V_{T_0})\in(\ell+j,\ell+j+1]\}]\leq c_0e^{-j}A_\theta(\ell+j).
    \end{align*}
    Now consider $B_\theta(\ell):=\mathbb{P}[\mathcal{E}_\theta\cap\log P_{T-T_0}^\pi f(W_{T_0})\in(\ell,\ell+1]]$. The $j=0$ case of the above equation then gives
    \begin{align*}
        B_\theta(\ell)\leq c_0A_\theta(\ell)+\mathbb{P}[\mathcal{E}_\theta\cap \log P_{T-T_0}^\pi f(V_{T_0})>\ell+1\geq\log P_{T-T_0}^\pi f(W_{T_0})].
    \end{align*}
    Now consider $D_{\mathcal{E}_\theta}$. It will satisfy (using $B=\{x\in G\mid \log P_{T-T_0}^\pi f(x)\in(\ell+1,\infty)\}$)
    \begin{align*}
        D_{\mathcal{E}_\theta}\geq|\mathbb{P}[\mathcal{E}_\theta\cap\{\log P_{T-T_0}^\pi f(W_{T_0})>\ell+1\}]-\mathbb{P}[\mathcal{E}_\theta\cap\{\log P_{T-T_0}^\pi f(V_{T_0})>\ell+1\}]|.
    \end{align*}
    We note that $\mathbb{P}[\mathcal{E}_\theta\cap\{\log P_{T-T_0}^\pi f(V_{T_0})>\ell+1\}]-\mathbb{P}[\mathcal{E}_\theta\cap\{\log P_{T-T_0}^\pi f(W_{T_0})>\ell+1\}]=\mathbb{P}[\mathcal{E}_\theta\cap \{\log P_{T-T_0}^\pi f(V_{T_0})>\ell+1\geq\log P_{T-T_0}^\pi f(W_{T_0})\}]-\mathbb{P}[\mathcal{E}_\theta\cap \{\log P_{T-T_0}^\pi f(W_{T_0})>\ell+1\geq\log P_{T-T_0}^\pi f(V_{T_0})\}]$, so we have the bound
    \begin{align*}
        &\mathbb{P}[\mathcal{E}_\theta\cap \{\log P_{T-T_0}^\pi f(V_{T_0})>\ell+1\geq\log P_{T-T_0}^\pi f(W_{T_0})\}]\\&\leq\mathbb{P}[\mathcal{E}_\theta\cap \{\log P_{T-T_0}^\pi f(W_{T_0})>\ell+1\geq\log P_{T-T_0}^\pi f(V_{T_0})\}]+D_{\mathcal{E}_\theta}.
    \end{align*}
    Thus
    \begin{align*}
        B_\theta(\ell)\leq c_0A_\theta(\ell)+\mathbb{P}[\mathcal{E}_\theta\cap \{\log P_{T-T_0}^\pi f(W_{T_0})>\ell+1\geq\log P_{T-T_0}^\pi f(V_{T_0})\}]+D_{\mathcal{E}_\theta}.
    \end{align*}
    We now partition the remaining probability
    \begin{align*}
        &\mathbb{P}[\mathcal{E}_\theta\cap \{\log P_{T-T_0}^\pi f(W_{T_0})>\ell+1\geq\log P_{T-T_0}^\pi f(V_{T_0})\}]\\
        &=\sum\limits_{j\geq 1}\mathbb{P}[\mathcal{E}_\theta\cap \{\ell+j<\log P_{T-T_0}^\pi f(W_{T_0})\leq\ell+j+1\}\cap\{\log P_{T-T_0}^\pi f(V_{T_0})\leq\ell+1\}]\\
        &\leq\sum\limits_{j\geq 1} c_0e^{-j}A_\theta(\ell+j).
    \end{align*}
    Therefore, $B_\theta(\ell)\leq c_0A_\theta(\ell)+D_{\mathcal{E}_\theta}+\sum\limits_{j\geq 1} c_0e^{-j}A_\theta(\ell+j)$. Finally, by using $B=\{x\in G\mid \log P_{T-T_0}^\pi f(x)\in(\ell,\ell+1]\}$ in the definition for $D$, we find $A_\theta(\ell)\leq  B_\theta(\ell)+D_{\mathcal{E}_\theta}$. This gives the desired bound:
    \begin{align*}
        (1-c_0)A_\theta(\ell)\leq c_0\sum_{j\geq1}e^{-j}A_\theta(\ell+j)+2D_{\mathcal{E}_\theta}.
    \end{align*}
\end{proof}

\section{Stopped score bound}\label{sec7}
In this section we prove Lemma~\ref{lem:stopped-score}. Our proof is very similar to the proof of \cite[Lemma 3.4]{lu2026weaktype}.
\begin{proof}[Proof of Lemma~\ref{lem:stopped-score}]
    For shorthand, let $\Phi(r)=r\log r-r+1$. By the jump It\^o formula,
    \begin{align*}
        &\log p_{T-t}(V_t)-\log p_{T-\theta}(V_\theta)-\int\limits_\theta^t\sum_{i=1}^N\sum_{g\in G_0}\pi(g)\Phi(\rho_{s,i,g}(V_{s^-}))\ ds
    \end{align*}
    is a martingale. By optional stopping,
    \begin{align*}
        \mathbb{E}[\log p_{T-\sigma}(V_\sigma)-\log p_{T-\theta}(V_\theta)\mid \widetilde{\mathcal{F}}_\theta]=\mathbb{E}\left[\int\limits_\theta^\sigma\sum\limits_{i=1}^N\sum\limits_{g\in G_0}\pi(g)\Phi(\rho_{s,i,g}(V_{s^-}))\ ds\Big|\ \widetilde{\mathcal{F}}_\theta\right].
    \end{align*}
    Since we have the $\rho$ bound $\kappa_{\tau,\pi}^{-1}\leq\rho_{t,i,g}(x)\leq\kappa_{\tau,\pi}$ by Lemma~\ref{lem:rho_bound} for $\tau=T-T_0$, we have a bound $(\rho_{t,i,g}(V_{t^-})-1)^2\leq\frac{2(\kappa_{\tau,\pi}-1)}{\log\kappa_{\tau,\pi}}\Phi(\rho_{t,i,g}(V_{t^-}))$ by Chen \cite[Equation 40]{chen2025talagrand}. Thus
    \begin{align*}
        &\mathbb{E}\left[\int\limits_\theta^\sigma\sum\limits_{i=1}^N\sum\limits_{g\in G_0}\pi(g)(\rho_{t,i,g}(V_{t^-})-1)^2\ dt\Big|\ \widetilde{\mathcal{F}}_\theta\right]\\&\qquad\leq\frac{2(\kappa_{\tau,\pi}-1)}{\log \kappa_{\tau,\pi}}\mathbb{E}[\log p_{T-\sigma}(V_\sigma)-\log p_{T-\theta}(V_\theta)\mid \widetilde{\mathcal{F}}_\theta].
    \end{align*}
    For the remaining expectation, if $\log P_{T-\theta}^\pi f(V_\theta)\geq\ell+1$, by definition of $\sigma$, $\sigma=\theta$, so we have zero. If instead $\ell\leq\log P_{T-\theta}^\pi f(V_\theta)<\ell+1$, we have $R_\theta=0$ and the difference is at most $1+\log\kappa_{\tau,\pi}$. Finally, if $\log P_{T-\theta}^\pi f(V_\theta)<\ell$, $\log P_{T-\theta}^\pi f(V_\theta)=\ell-R_\theta$ and the difference is at most $1+R_\theta+\log \kappa_{\tau,\pi}$. Therefore
    \begin{align*}
        \mathbb{E}\left[\int\limits_\theta^\sigma\sum\limits_{i=1}^N\sum\limits_{g\in G_0}\pi(g)(\rho_{t,i,g}(V_{t^-})-1)^2\ dt\Big|\ \widetilde{\mathcal{F}}_\theta\right]\leq\frac{2(\kappa_{\tau,\pi}-1)}{\log\kappa_{\tau,\pi}}(1+R_\theta+\log\kappa_{\tau,\pi}).
    \end{align*}
    For the second statement, note that on $\mathcal{E}_\theta$, $\bar{\delta}=\frac{\alpha}{R_\theta+1}$, so by the above bound,
    \begin{align*}
        \mathscr{S}_A&\leq\frac{2\alpha^2\kappa_{\tau,\pi}}{\Lambda_{\tau,\pi}}\mathbb{E}\left[\frac{\mathbbold{1}_A}{R_\theta+1}\right]+2\alpha^2(\kappa_{\tau,\pi}-1)\mathbb{E}\left[\frac{\mathbbold{1}_A}{(R_\theta+1)^2}\right].
    \end{align*}
\end{proof}

\section{Additional Proofs}
\begin{lemma}\label{generator}
    The generator of the semigroup $(P_t^\pi)_{t\geq0}$ is
    \begin{align*}
        Lf(x)=\sum\limits_{i=1}^N\sum\limits_{g\in G_0}\pi(g)(f(x^{i\to g})-f(x)).
    \end{align*}
\end{lemma}
\begin{proof}
    For each coordinate $i$, define $E_i u(x):=\sum\limits_{g\in G_0}\pi(g)u(x^{i\to g})$. Independence of the coordinate refreshes gives
    \begin{align*}
        P_t^\pi=\prod_{i=1}^N\bigl(e^{-t}I+(1-e^{-t})E_i\bigr).
    \end{align*}
    Differentiating at $t=0$, we obtain
    \begin{align*}
        Lu(x)=\sum_{i=1}^N(E_i-I)u(x)=\sum_{i=1}^N\sum_{g\in G_0}\pi(g)\bigl(u(x^{i\to g})-u(x)\bigr).
    \end{align*}
\end{proof}

\begin{lemma}\label{self_adjoint}
    The generator $L$ of the semigroup $(P_t^\pi)_{t\geq0}$ is self-adjoint on $L^2(\pi^{\otimes N})$.
\end{lemma}
\begin{proof}
    Under $\mu=\pi^{\otimes N}$, the operator $E_i$ from Lemma~\ref{generator} is conditional expectation given all coordinates except the $i$th. Hence $E_i$ is an orthogonal projection on $L^2(\mu)$ and is self-adjoint. Therefore $L=\sum\limits_{i=1}^N(E_i-I)$ is self-adjoint.
\end{proof}

\begin{lemma}\label{reverse_generator}
    The time-dependent generator of the time-reversed version of the refresh process starting at $f\pi^{\otimes N}$, $0\leq t\leq T$, $f:G\to(0,\infty)$ with $\int_G f\ d\pi^{\otimes N}=1$, for any $T>0$, is
    \begin{align*}
        \widetilde{L}_tu(x)=\sum\limits_{i=1}^N\sum\limits_{g\in G_0}\pi(g)\left(\frac{p_{T-t}(x^{i\to g})}{p_{T-t}(x)}\right)(u(x^{i\to g})-u(x)),
    \end{align*}
    where $p_t=P_t^\pi f$ is the density of the forward process's law with respect to $\pi^{\otimes N}$.
\end{lemma}
\begin{proof}
    Write $\mu=\pi^{\otimes N}$ and $p_s=P_s^\pi f$, so the forward process has law $p_s\mu$. If $q(x,y)$ denotes its jump rate, Bayes' formula gives the reverse rate
    \begin{align*}
        \widetilde q_t(x,y)=\frac{p_{T-t}(y)\mu(y)}{p_{T-t}(x)\mu(x)}q(y,x),\qquad x\neq y.
    \end{align*}
    For $y=x^{i\to g}\ne x$, we have $q(y,x)=\pi(x_i)$ and $\mu(y)/\mu(x)=\pi(g)/\pi(x_i)$. Thus
    \begin{align*}
        \widetilde q_t(x,x^{i\to g})=\pi(g)\frac{p_{T-t}(x^{i\to g})}{p_{T-t}(x)}.
    \end{align*}
    This gives generator
    \begin{align*}
        \widetilde L_tu(x)=\sum_{i=1}^N\sum_{g\in G_0}\pi(g)\frac{p_{T-t}(x^{i\to g})}{p_{T-t}(x)}\bigl(u(x^{i\to g})-u(x)\bigr).
    \end{align*}
\end{proof}

\begin{lemma}\label{lem:rho_bound}
    For all $f:G\to(0,\infty)$ with $\int_Gf\ d\pi^{\otimes N}=1$, $\tau>0$, $0\leq t\leq T-\tau$, $x\in G$, $i\in\{1,\cdots,N\}$, $g\in G_0$,
    \begin{align*}
        \kappa_{\tau,\pi}^{-1}\leq\rho_{t,i,g}(x)\leq\kappa_{\tau,\pi},
    \end{align*}
    where $\kappa_{t,\pi}=\frac{e^{-t}+(1-e^{-t})\max\limits_{y\in G_0}\pi(y)}{(1-e^{-t})\min\limits_{y\in G_0}\pi(y)}$ and $\rho_{t,i,g}(x)=\frac{p_{T-t}(x^{i\to g})}{p_{T-t}(x)}$, $p_t(x)=P_t^\pi f(x)$.
\end{lemma}
\begin{proof}
    By definition of $P_t^\pi$,
    \begin{align*}
        \rho_{t,i,g}(x)=\frac{P_{T-t}^\pi f(x^{i\to g})}{P_{T-t}^\pi f(x)}=\frac{\int_G f(y)\ d\pi_{T-t,x^{i\to g}}(y)}{\int_G f(y)\ d\pi_{T-t,x}(y)}=\frac{\sum\limits_{y\in G}f(y)\pi_{T-t,x^{i\to g}}(y)}{\sum\limits_{y\in G}f(y)\pi_{T-t,x}(y)}.
    \end{align*}
    Since $f>0$, this can be bounded by the finite quantity (since $\pi$ has full support)
    \begin{align*}
        \rho_{t,i,g}(x)\leq\max\limits_{y\in G}\frac{\pi_{T-t,x^{i\to g}}(y)}{\pi_{T-t,x}(y)}.
    \end{align*}
    Then noting that $x^{i\to g}$ and $x$ only differ in coordinate $i$ and the measures match on all other coordinates,
    \begin{align*}
        &\rho_{t,i,g}(x)\leq\max\limits_{y\in G}\frac{e^{-(T-t)}\mathbbold{1}_{g=y_i}+(1-e^{-(T-t)})\pi(y_i)}{e^{-(T-t)}\mathbbold{1}_{x_i=y_i}+(1-e^{-(T-t)})\pi(y_i)}\leq\max_{a\in G_0}\frac{e^{-(T-t)}+(1-e^{-(T-t)})\pi(a)}{(1-e^{-(T-t)})\pi(a)}\\&\leq\frac{e^{-(T-t)}+(1-e^{-(T-t)})\max\limits_{y\in G_0}\pi(y)}{(1-e^{-(T-t)})\min\limits_{y\in G_0}\pi(y)}\leq\kappa_{\tau,\pi}.
    \end{align*}
    The lower bound follows by noting $\frac{1}{\rho_{t,i,g}(x)}=\rho_{t,i,x_i}(x^{i\to g})$.
\end{proof}

\begin{lemma}\label{lem:rep}
    Let $\varphi:G\to[0,1]$ and $\mu=\bigotimes\limits_{i=1}^N\mu_i$ be a product measure on $G$ factoring over the $G_0$ components. For $X\sim \mu$, define $\bar{\varphi}:=\mathbb{E}^\mu[\varphi(X)]$ and, for every $u\in G_0$, define
    \begin{align*}
        \varphi_i(u):=\int\varphi(x_1,\ldots,x_{i-1},u,x_{i+1},\ldots,x_N)\prod_{j\neq i}\mu_j(dx_j).
    \end{align*}
    Then
    \begin{align*}
        \sum\limits_{i=1}^N\sum\limits_{u\in G_0}\mu_i(u)(\varphi_i(u)-\bar{\varphi})^2\leq\mathrm{Var}_\mu(\varphi)\leq\frac{1}{4}.
    \end{align*}
\end{lemma}
\begin{proof}
    Let $X\sim\mu$ and set $Z:=\varphi(X)-\bar\varphi$,$g_i:=\varphi_i(X_i)-\bar\varphi=\mathbb{E}[Z\mid X_i]$. The $g_i$ are centered and depend on independent coordinates, so $\mathbb{E}[g_i g_j]=0$ for $i\neq j$. Moreover, $\mathbb{E}[Zg_i]=\mathbb{E}\bigl[g_i\mathbb{E}[Z\mid X_i]\bigr]=\mathbb{E}[g_i^2]$. Consequently,
\begin{align*}
    0\leq\mathbb{E}\left[\left(Z-\sum_{i=1}^N g_i\right)^2\right]=\operatorname{Var}_\mu(\varphi)-\sum_{i=1}^N\mathbb{E}[g_i^2].
\end{align*}
Since $\mathbb{E}[g_i^2]=\sum\limits_{u\in G_0}\mu_i(u)\bigl(\varphi_i(u)-\bar\varphi\bigr)^2$, this proves the first inequality. Finally, $0\le\varphi\le1$ implies $\operatorname{Var}_\mu(\varphi)\leq\bar\varphi-\bar\varphi^2\leq\frac14$.
\end{proof}

\section{Statement on AI Use}
Artificial Intelligence - specifically ChatGPT-6.1 Sol - was used to formalize the proof in Lean 4. The formalization was then verified by the author to confirm that the correct theorem was proven. AI was also used for language proofreading and to improve exposition and clarity.

\providecommand{\bysame}{\leavevmode\hbox to3em{\hrulefill}\thinspace}

\end{document}